\documentclass[preprint,12pt]{elsarticle}

\usepackage{graphicx}
\usepackage{amsfonts}
\usepackage[T1]{fontenc}
\usepackage{amsmath}
\usepackage{amssymb}
\usepackage{bm}
\usepackage{optidef}
\usepackage{hyperref}
\usepackage{csquotes}
\usepackage{bm}
\usepackage{algorithm}
\usepackage{algpseudocode}
\usepackage{circuitikz}
\usepackage{siunitx}
\usepackage{booktabs}
\usepackage{subfiles}
\usepackage{subcaption}
\usepackage{tikz,tikz-3dplot}
\usepackage{xcolor}

\DeclareSIUnit{\rad}{rad}
\let\vec\bm

\let\mat\mathbf
\def\param{\vec{\chi}}

\newcommand\figref{Figure~\ref}

\let\amssymbboxplus\boxplus
\let\amssymbboxminus\boxminus

\renewcommand{\boxplus}{\mathbin{\mathop\amssymbboxplus}}
\renewcommand{\boxminus}{\mathbin{\mathop\amssymbboxminus}}

\newcommand{\norm}[1]{\left\lVert#1\right\rVert}

\renewcommand{\algorithmicrequire}{\textbf{Input:\,}}
\DeclareMathOperator*{\argmax}{arg\,max}
\DeclareMathOperator*{\diag}{diag}

\numberwithin{equation}{section}

\def\linvel{\mathrm{v}}
\def\linacc{\dot{\vec{v}}}
\def\angvel{\vec{\omega}}
\def\angacc{\dot{\vec{\omega}}}
\def\rv{\psi}
\def\quat{e}
\def\conf{q}
\def\jaca{\mat{A}}
\def\jacb{\mat{B}}
\def\tA{\text{A}}
\def\tB{\text{B}}
\def\tP{\text{P}}
\def\tc{\text{c}}
\def\tg{\text{g}}
\def\th{\text{h}}
\def\step{\phi}
\def\lambdah{\vec{\lambda}_{\text{h}}}
\def\lambdanh{\vec{\lambda}_{\text{nh}}}
\def\lambdan{\vec{\lambda}_{\text{n}}}
\def\lambdat{\vec{\lambda}_{\text{t}}}
\def\lambdaa{\vec{\lambda}_{\text{A}}}
\def\lambdab{\vec{\lambda}_{\text{B}}}
\def\drydiagonal{\mat{\Sigma}_\text{B}}

\def\combresiduals{ \Delta \vec{z}}
\def\stateweight{\kappa}
\def\statenoise{\eta}
\def\obsnoise{\xi}

\makeatletter
\renewcommand*\env@matrix[1][*\c@MaxMatrixCols c]{%
    \hskip -\arraycolsep
    \let\@ifnextchar\new@ifnextchar
    \array{#1}}
\makeatother

\newcommand{\complementarity}[2]{#1^T #2 = \vec{0}, \qquad #1, #2 \geq \vec{0}}
\newcommand{\complementarityalign}[2]{{#1}^T{#2}=\vec{0}, &\qquad #1, #2 \geq \vec{0}}

\begin{document}

\begin{frontmatter}
    \title{Joint parameter and state estimation for regularized time-discrete multibody dynamics}
    
    \author[]{Hannes Marklund}
    \author[]{Martin Servin}
    \author[]{Mats G. Larson}
    
    \affiliation[label1]{
                city={Umeå},
                postcode={SE-90187}, 
                country={Sweden}}
    \affiliation[org1]{organization={Department of Physics, Umeå University}}
    \affiliation[org2]{organization={Department of Mathematics and Mathematical Statistics, Umeå University}}

    \begin{abstract}
        We develop a method for offline parameter estimation of discrete multibody dynamics with regularized and frictional kinematic constraints.
        This setting leads to unobserved degrees of freedom, which we handle using joint state and parameter estimation.
        Our method finds the states and parameters as the solution to a nonlinear least squares optimization problem based on the inverse dynamics and the observation error.
        The solution is found using a Levenberg--Marquardt algorithm with derivatives from automatic differentiation and custom differentiation rules for the complementary conditions that appear due to dry frictional constraints.
        
        We reduce the number of method parameters to the choice of the time-step, regularization coefficients, and a parameter that controls the relative weighting of inverse dynamics and observation errors.
        We evaluate the method using synthetic and real measured data, focusing on performance and sensitivity to method parameters.
        In particular, we optimize over a 13-dimensional parameter space, including inertial, frictional, tilt, and motor parameters, using data from a real Furuta pendulum.
        Results show fast convergence, in the order of seconds, and good agreement for different time-series of recorded data over multiple method parameter choices.
        However, very stiff constraints may cause difficulties in solving the optimization problem.
        We conclude that our method can be very fast and has method parameters that are robust and easy to set in the tested scenarios.
    \end{abstract}
    

    
    
    \end{frontmatter}

\section{Introduction}
\label{sec:introduction}
Simulators are an important tool for the design and control of robots, vehicles, and other cyber-physical systems.
Often, the simulators are required to run in real-time or faster to support hardware-in-the-loop, operator training, model-based control, or reinforcement learning.
To maximize the computational speed, these simulators are based on coarse-grained models that resolve the dynamics on the time scale at which the controllers operate, but not finer.
This usually amounts to using multibody system dynamics with joints and contacts represented as kinematic constraints and complementarity conditions \cite{baraff1996linear,stewart1996implicit} integrated at fixed time-step around 10-\SI{1000}{Hz}.
For numerically robust incorporation of kinematic loops and frictional contacts, it is customary to explicitly represent all six degrees of freedom per rigid body (maximal coordinates), and solve for regularized constraint forces using Lagrange multipliers \cite{lacoursiere2007ghosts}.
A common hurdle to using simulations is that they often involve parameters that are not derivable from first-principle models or identifiable by direct measurements.
Instead, these parameters must be determined from observations of the full-system dynamics \cite{lee2023robot}.
If task-relevant parameters are not accurately determined, the result of simulations will not transfer to the real world.
This so-called reality-gap is a major obstacle to using simulation to control robots, directly or using AI-models pre-trained on simulated data \cite{choi2021use}.

In parameter identification, we seek model parameters that minimize the discrepancy between the model and observation data from the real system, resulting in methods that are often divided into two groups: The \textit{equation error method} and the \textit{simulation error method} \cite{lee2023robot}.
In the simulation error method, we compare real and simulated time-series of observables.
This approach typically results in good global behavior at the expense of a highly non-convex optimization landscape.
Without some stabilizing mechanism, such as closed-loop control, a satisfactory solution might not exist.
In the equation error method, we compare local predictions using the dynamic equations.
In discrete mechanics, we can predict the next state (forward dynamics) or predict the force required to follow a prescribed state transition (inverse dynamics).
The advantage is a more tractable optimization problem at the expense of demanding full observability---the property that observables map directly to the state variables in the equations.

In regularized time-discrete multibody dynamics, we allow deviations from the otherwise ideal kinematic constraints to gain numerical robustness.
There are two reasons for constraint deviations:
Local linearization of the constraint equations and regularization, which may be associated with actual physical properties.
These deviations pose a problem for observability--constrained degrees of freedom are typically not measured.
Furthermore, even if the constrained degrees of freedom are measured, the deviations may be a result of numerical artifacts that do not correspond to the physical deviation.
Since the simulation error method works for partially observable systems, it remains the standard choice.
However, problems may arise if we put constraints on optimization time or if the dynamics is unstable.
In such cases, using the equation error method would be desirable.
To address the observability problem, we need to estimate the states alongside the parameters.
State estimation of multibody dynamics has been investigated extensively \cite{nayaKalmanFiltersBased2023,pastorino2013state, sanjurjoAccuracyEfficiencyComparison2017} and combined estimation to a lesser extent \cite{nayaKalmanFiltersBased2023, blanchardPolynomialChaosBasedKalman2010}.
However, existing approaches consider online estimation, meaning time-varying parameters. This is not a suitable fit for offline estimation---we do not want the parameters to compensate for model errors locally.

In this work, we develop an offline joint state and parameter estimation method for regularized time-discrete multibody dynamics with frictional joints.
We find the states and the parameters by solving a nonlinear least squares optimization problem based on the inverse dynamics and observation error.
We study and evaluate the computational efficiency and sensitivity of the method in examples with both synthetic and real measured data.
In particular, we record measurements using a real Furuta pendulum experimental setup.
To compute derivatives, we use our own differentiable simulator with custom differentiation rules for the complementarity conditions that appear when dry friction is included. 
The simulator is implemented using the Python library \textsc{Jax} \cite{jax2018github}.

The remainder of the paper is organized as follows:
Section \ref{sec:multibody_dynamics} presents the relevant background regarding regularized time-discrete multibody dynamics.
Section \ref{sec:parameter_identification} covers our combined parameter and state estimation method.
In Section \ref{sec:case_study}, we present the experimental and numerical studies, including results.
We summarize our findings, draw conclusions, and outline future research directions in Section \ref{sec:conclusions}.

\section{Time-discrete multibody dynamics}
\label{sec:multibody_dynamics}
This section covers multibody dynamics with kinematic constraints and frictional contacts.
We present a time-stepper based on a maximal coordinate description and regularized constraints to allow numerically robust incorporation of kinematic loops and frictional contacts.
We then describe our constraint solver and how we use time-discrete multibody dynamics in optimization.

\subsection{Representation}
We are concerned with systems well-approximated by rigid bodies interacting with joints and contacts.
The position of all points in a rigid body can be determined from its shape, the position of a reference point $\vec{r} \in \mathbb{R}^3$, and an orientation, which we choose to represent by a unit quaternion $\vec{\quat} \in \mathbb{H}_0$ from the body frame to the global inertial frame.
We prefer writing $\mathbb{H}_0$ over $\mathbb{R}^4$ to emphasize the additional structure and allow quaternion multiplication, which we denote with $\star$.
We write the \emph{configuration} of a rigid body $\tA$ as the tuple
\begin{equation}
    \label{eq:configuration}
    \vec{\conf}_\tA= (\vec{r}_\tA, \vec{\quat}_\tA) \in \mathbb{R}^{3} \times \mathbb{H}_0,
\end{equation}
where $\vec{r}_\tA$ points to $\tA$'s center of mass.
The generalized velocity of $\tA$ is given by
\begin{equation}
    \label{eq:generalized_velocity}
    \vec{v}_A= \begin{bmatrix}
        \vec{\linvel}_A^T & \angvel_A^T
    \end{bmatrix}^T\in \mathbb{R}^{6},
\end{equation}
where $\vec{\linvel}_\tA$ is the linear velocity and  $\angvel_\tA$ is the rotational velocity.
To describe the time evolution of $\tA$, we use the Newton-Euler equations
\begin{equation}
    \underbrace{
        \begin{bmatrix}
            \vec{\mathrm{f}}_\tA \\ \vec{\tau}_\tA
        \end{bmatrix}
    }_{\vec{f}_\tA} =
    \underbrace{
        \begin{bmatrix}
            m_\tA \mat{I}_{3\times3} & \vec{0}   \\
            \vec{0}     & \mat{J}_\tA(\vec{\conf}_\tA)
        \end{bmatrix}
    }_{\mat{M}_\tA(\vec{\conf}_\tA)}
    \begin{bmatrix}
        \linacc_\tA \\ \angacc_\tA
    \end{bmatrix}
    +
    \underbrace{
        \begin{bmatrix}
            \vec{0} \\
            \angvel_\tA \times \mat{J}_\tA(\vec{\conf}_\tA) \angvel_\tA
        \end{bmatrix}
    }_{-\vec{f}_{\tg,\tA}},
\end{equation}
where we use the following notation:
First, the mass $m_\tA$ times the identity matrix $\mat{I}_{3\times3}$ and configuration-dependent inertia matrix $\mat{J}_\tA$ form the mass matrix $\mat{M}_\tA$.
Second, the force $\vec{\mathrm{f}}_\tA$ and torque $\vec{\tau}_\tA$ forms the generalized force $\vec{f}_\tA$.
Finally, $\vec{f}_{\tg,\tA}$ denotes the generalized fictitious force.
We use dot-notation to denote derivatives with respect to time.
To compactly describe a system of multiple bodies, we drop the indices
\begin{equation}
    \label{eq:newton_euler_multiple_bodies}
    \underbrace{
        \begin{bmatrix}
            \vec{f}_\tA \\ \vec{f}_\tB \\ \vdots
        \end{bmatrix}
    }_{\vec{f}}
    =
    \underbrace{
        \begin{bmatrix}
            \mat{M}_\tA(\vec{\conf}_\tA) & \mat{0}   & \hdots \\
            \mat{0}   & \mat{M}_\tB (\vec{\conf}_\tB) & \hdots \\
            \vdots    & \vdots    & \ddots
        \end{bmatrix}
    }_{\mat{M}(\vec{\conf})}
    \underbrace{
        \begin{bmatrix}
            \dot{\vec{v}}_\tA \\ \dot{\vec{v}}_B \\ \vdots
        \end{bmatrix}
    }_{\dot{\vec{v}}}
    +
    \underbrace{
        \begin{bmatrix}
            \vec{f}_{\tg, \tA} \\ \vec{f}_{\tg, \tB}  \\ \vdots
        \end{bmatrix}
    }_{-\vec{f}_g}.
\end{equation}
The configuration $\vec{\conf}$ of a system of $n_b$ rigid bodies resides in $\mathcal{M} = (\mathbb{R}^3 \times \mathbb{H}_0)^{n_b}$.

\subsection{Kinematic Constraints}
We describe joints and contacts between bodies as kinematic constraints---constraints that restrict the movement of the associated bodies.
The most basic kinematic constraint is a \textit{holonomic constraint}.
A holonomic constraint, or set of holonomic constraints, takes the form $\vec{g}(\vec{\conf}) = \vec{0}$, where $\vec{g}: \mathcal{M} \rightarrow \mathbb{R}^{d_{\th}}$ is a continuously differentiable function.
For a 3D multibody system, it is convenient to split $\vec{g}$ into $n_h$ elementary constraints $\vec{g} = \begin{bmatrix} \vec{g}_1^T \hdots \vec{g}_{n_h}^T \end{bmatrix}^T$, where the $i$th constraint maps $\mathcal{M}$ to $\mathbb{R}^{d_{\th, i}}$ and $d_{\th, i}\in\{1,\hdots,6\}$.
We write the equations of motion as
\begin{align}
    \label{eq:manipulator}
    \mat{M} (\vec{\conf}) \dot{\vec{v}} & = \vec{f} + \vec{f}_g + \mat{G}(\vec{\conf})^T \vec{\lambda}, \\
    \label{eq:manipulator2}
    \vec{g}(\vec{\conf})                & = \vec{0},
\end{align}
a differential-algebraic system of equations (DAE).
Here, $\mat{G}: \mathcal{M} \to \mathbb{R}^{d \times 6n_b}$ is the Jacobian matrix of $\vec{g}$, defined by $\dot{\vec{g}}(\vec{\conf}) = \mat{G}(\vec{\conf}) \vec{v}$, and $\vec{\lambda}$ are Lagrange multipliers, which couples to the constraint force $\mat{G}(\vec{\conf})^T \vec{\lambda}$ necessary for maintaining the constraint.

Other kinds of kinematic constraints are also useful.
Of relevance here are two other types, namely \emph{nonholonomic constraints} and \emph{contact constraints}. Nonholonomic constraints are constraints on the form $\bar{\vec{g}}(\vec{\conf}, \vec{v}, t) = \vec{0}$ that one cannot reduce to a holonomic constraint.
However, it is enough for us only to consider nonholonomic constraints that are linear in velocity, that is $\bar{\mat{G}}(\vec{\conf}) \vec{v} = \vec{\gamma}(t)$.
These are useful for modeling friction after including regularization.
Contact constraints describe non-penetration conditions and take the form $\vec{c}(\vec{\conf}) \geq \vec{0}$.

\subsection{Time-stepping}
In the presence of holonomic constraints, the simulator has to integrate the DAE of \eqref{eq:manipulator} and \eqref{eq:manipulator2}.
As a starting point, one may differentiate the algebraic constraints to convert the DAE into an ODE.
Equations \eqref{eq:manipulator} and \eqref{eq:manipulator2} imply
\begin{equation}
    \label{eq:dae_to_ode}
    \begin{bmatrix}
        \mat{M} & -\mat{G}^T \\
        \mat{G} & \vec{0}    \\
    \end{bmatrix}
    \begin{bmatrix}
        \dot{\vec{v}} \\ \vec{\lambda}
    \end{bmatrix}
    =
    \begin{bmatrix}
        \vec{f} + \vec{f}_g \\
        -\dot{\mat{G}}\vec{v}
    \end{bmatrix},
\end{equation}
where the holonomic constraint $\vec{g}(\vec{\conf}) = \vec{0}$ has been replaced by the acceleration constraint $\ddot{\vec{g}}(\vec{\conf}) = \mat{G}\dot{\vec{v}} + \dot{\mat{G}}\vec{v} = \vec{0}$.
This is exact for initial positions and velocity in compliance with the constraint.
However, when integrating numerically, errors will cause bodies to drift away from the constraint surface $\vec{g}(\vec{q}) = \vec{0}$.
A common solution to this problem is to include so-called Baumgarte stabilization \cite{baumgarte1972stabilization}, which makes the holonomic constraints function as damped springs.
Another problem is that the linear system can be ill-conditioned or even singular due to large mass ratios, closed loops, or redundant constraints.
To prevent ill-conditioning, we can add a diagonal perturbation (regularization) to the constraint equations, which also softens the constraints.

In this work, we use an integration scheme with fixed time-step introduced by Lacoursiere, known as \textsc{Spook}  \cite{lacoursiere2007ghosts}.
A fixed time-step is desirable for the predictable performance required in real-time applications. \textsc{Spook} is a semi-implicit extension of the Verlet stepper with time-symmetric linearization of the constraint forces.
Verlet and \textsc{Spook} are both \emph{variational integrators}, meaning they preserve phase-space volume and have good energy and momentum-preserving behavior \cite{hairer2006geometric}.Variational integrators produce \emph{shadowing trajectories}, trajectories that remain close to the trajectory of the continuous dynamics, limited by chaotic dynamics and phase space saddle points \cite{chandramoorthy2021probability}.
\textsc{Spook} uses a combination of Baumgarte stabilization and regularization that can be derived using discrete variational mechanics.
Each elementary holonomic constraint $\vec{g}_i$ is now associated with a compliance scalar $\epsilon_{i}$ and a damping scalar $\tau_{i}$. The compliance $\epsilon_{i}$ can be interpreted as an inverse spring coefficient and the damping $\tau_{i}$ is a relaxation time.
We can either tune these parameters for optimal constraint satisfaction, to improve conditioning, or set them to represent the real physical parameters.
For a multibody system with $n_\th$ holonomic constraints, the \textsc{Spook} stepper reads
\begin{equation}
    \label{eq:spook_holonomic}
    \begin{bmatrix}
        \mat{M}_k & -\mat{G}_k^T \\
        \mat{G}_k & \vec{\Sigma}_\th
    \end{bmatrix}
    \begin{bmatrix}
        \vec{v}_{k+1} \\
        \lambdah
    \end{bmatrix}
    =
    \begin{bmatrix}
        \mat{M}_k \vec{v}_k - h (\vec{f}_k + \vec{f}_{\tg,k}) \\
        -\frac{4}{h} \vec{\Upsilon}_\th \vec{g}_k + \vec{\Upsilon}_\th \mat{G}_k \vec{v}_k
    \end{bmatrix}, \\
\end{equation}
where $h$ is the time-step, $\lambdah$ is the vector of multipliers associated with holonomic constraints, and
\begin{align}
    \label{eq:upsilon}
    \gamma_{i} & = \frac{1}{1 + \frac{4 \tau_{i}}{h}}, \qquad i=1,\hdots,n_\th,            \\
    \mat{\Sigma}_\th   & =  \frac{4}{h^2}\text{diag}(
    \epsilon_{1} \gamma_{1} \mat{I}_{d_1 \times d_1},
    \epsilon_{2} \gamma_{2} \mat{I}_{d_2 \times d_2},
    \hdots,
    \epsilon_{n_\th} \gamma_{n_\th} \mat{I}_{d_{n_h} \times d_{n_h}}
    ),                                                                                           \\
    \mat{\Upsilon}_\th & = \text{diag}(\gamma_{1} \mat{I}_{d_1 \times d_1}, \gamma_{2} \mat{I}_{d_2 \times d_2}, \hdots,\gamma_{n_\th} \mat{I}_{d_{n_h} \times d_{n_h}}).
\end{align}
The subscript $k$ is used to denote discrete time, $t=hk$, and we use the short-hand notataion $\vec{g}_k = \vec{g}(\vec{\conf}_k), \mat{G}_k = \mat{G}(\vec{\conf}_k), \mat{M}_k = \mat{M}(\vec{\conf}_k)$.
It is worth noting that $\lambdah$ has absorbed a factor $h$ and is now an impulse rather than a force (compare with $\vec{\lambda}$ in \eqref{eq:dae_to_ode}).

The new velocity $\vec{v}_{k+1}$ is used to update the configuration from $\vec{q}_k$ to $\vec{q}_{k+1}$.
The translation is updated by
\begin{equation}
    \vec{r}_{\tA, {k+1}} = \vec{r}_{\tA, k} + h\vec{v}_{\tA, k+1},
\end{equation}
and similarly for rigid bodies $\tB$, $\text{C}$, etc.
The rotation has a similar update rule but requires the introduction of new notation to state compactly.
Note that $h \vec{\omega}_{\tA, k+1}$ can be seen as the product of a unit vector and an angle, a so-called \emph{rotation vector}.
We want to rotate $\vec{e}_{\tA, k}$ around that unit vector by that angle.
For a general rotation vector $\vec{\rv}$ applied to a unit quaternion $\vec{\quat}$, the result $\vec{\quat}'$ is given by
\begin{equation}
    \vec{\quat}' = 
    \vec{\quat} \star
    \begin{bmatrix}
        \cos (\frac{\|\vec{\rv}\|}{2}) \\
        \frac{\vec{\rv}}{\|\vec{\rv}\|} \sin (\frac{\|\vec{\rv}\|}{2})
    \end{bmatrix},
\end{equation}
where $\star$ denotes quaternion multiplication.
This operation is common enough to warrant the introduction of the operators
\begin{eqnarray}
    \label{eq:boxdef}
    \vec{\quat} \boxplus \vec{\rv} &= \vec{\quat}', \\
    \vec{\quat}' \boxminus \vec{\quat} &= \vec{\rv}.
\end{eqnarray}
This is a common way to extend the plus and minus operators for Lie groups such as rotation matrices and unit quaternions \cite{sola2017quaternion}.
Further, we allow use on tuples of multiple rotations and translations where they collapse to standard addition or subtraction.
The configuration update is now compactly written as
\begin{equation}
    \vec{q}_{k+1} = \vec{q}_k \boxplus h\vec{v}_{k+1}.
\end{equation}
Nonholonomic constraints that are linear in velocity, $\bar{\mat{G}}(\conf) \vec{v} = \vec{\gamma}(t)$ can be included in \textsc{Spook} using the discretization
\begin{equation}
    \label{eq:nonholonomic_disc}
    \mat{\Sigma}_\text{nh} \lambdanh + \bar{\mat{G}}_k \vec{v}_{k+1}= \vec{\gamma}_k,
\end{equation}
where the matrix $\mat{\Sigma}_\text{nh}$ provides regularization and  $\lambdanh$ is the vector of multipliers associated with nonholonomic constraints.
Together with holonomic constraints, the velocity update is given by
\begin{equation}
    \label{eq:spook_wihtout_complementarity}
    \begin{bmatrix}
        \mat{M}_k & -\mat{G}_k^T & -\bar{\mat{G}}_k^T    \\
        \mat{G}_k   & \mat{\Sigma}_\th & \mat{0} \\
        \bar{\mat{G}}_k & \mat{0} & \mat{\Sigma}_\text{nh} \\
    \end{bmatrix}
    \begin{bmatrix}
        \vec{v}_{k+1} \\
        \vec{\lambda}_\th \\
        \vec{\lambda}_\text{nh} \\
    \end{bmatrix}
        =
    \begin{bmatrix}
        \mat{M}_k \vec{v}_k - h (\vec{f}_k + \vec{f}_{\tg,k}) \\
        -\frac{4}{h} \vec{\Upsilon} \vec{g}_k + \vec{\Upsilon} \mat{G}_k \vec{v}_k \\
        \vec{\gamma}_k
    \end{bmatrix}.
\end{equation}
For the sake of less cluttered notation, we write the system compactly as
\begin{equation}
    \label{eq:spook_wihtout_complementarity_compact}
    \begin{bmatrix}
        \mat{M}_k & -\jaca_k^T     \\
        \jaca_k   & \mat{\Sigma}_\tA
    \end{bmatrix}
    \begin{bmatrix}
        \vec{v}_{k+1} \\
        \vec{\lambda}_\tA
    \end{bmatrix}
    =
    -\begin{bmatrix}
        \vec{p}_k \\
        \vec{a}_k
    \end{bmatrix},
\end{equation}
where $\mat{\Sigma}_\tA$, $\mat{A}_k$, $\vec{\lambda}_\tA$, $\vec{p}_k$, and $\vec{a}_k$ are given by direct comparison with \eqref{eq:spook_wihtout_complementarity}.

\subsection{Complementarity}
\label{sec:complementarity}
Contact constraints can be stated algebraically as $\vec{c}(\vec{\conf}) \geq \vec{0}$ where $\vec{c}: \mathcal{M} \rightarrow \mathbb{R}^{d_{\text{c}}}$ is a \emph{gap function}---it tells the displacement and rotation allowed before contact is made.
The continuous equations of motion with contact constraints read
\begin{equation}
    \begin{aligned}
        \mat{M} (\vec{\conf}) \dot{\vec{v}} & = \vec{f} + \vec{f}_\tg + \mat{C}(\vec{\conf})^T \vec{\lambda}, \\
        \complementarityalign{\vec{\lambda}}{\vec{c}(\vec{\conf})},
    \end{aligned}
\end{equation}
where $\mat{C}$ is the Jacobian matrix of $\vec{c}$.
The conditions $\vec{\lambda}^T\vec{c}(\vec{\conf}) = \vec{0}$ and $\vec{\lambda}, \vec{c}(\vec{\conf}) \geq \vec{0}$ are collectively called \emph{complementarity conditions}.
These ensure that at most one of $\vec{\lambda}$ and $\vec{c}(\vec{\conf})$ are nonzero and positive component-wise.

Recall that the time-stepper considers regularized constraints.
That is, we must associate the $n_\tc$ elementary contact constraints with compliance and damping terms, $\epsilon_{\tc, i}$ and $\tau_{\tc, i}$.
The complementarity conditions carry over to the time-stepper after discretization:
\begin{equation}
    \begin{aligned}
        \label{eq:spook_inequality}
        \begin{bmatrix}
            \mat{M}_k & -\mat{C}_k^T \\
            \mat{C}_k & \vec{\Sigma}_\tc
        \end{bmatrix}
        \begin{bmatrix}
            \vec{v}_{k+1} \\
            \vec{\lambda}_{\tc}
        \end{bmatrix}
         & +
        \begin{bmatrix}
            -\mat{M}\vec{v}_k + h (\vec{f} + \vec{f}_g) \\
            \frac{4}{h} \mat{\Upsilon}_\tc \vec{c}_k - \mat{\Upsilon}_\tc \mat{C}_k \vec{v}_k
        \end{bmatrix}
        =
        \begin{bmatrix}
            \vec{0} \\ \vec{w}_\tc
        \end{bmatrix}
        ,    \\
        \complementarityalign{\vec{\lambda}_\tc}{\vec{w}_\tc}.
    \end{aligned}
\end{equation}
Here, $\mat{\Upsilon}_\tc$ and $\mat{\Sigma}_\tc$ are given by \eqref{eq:upsilon} and replacing holonomic damping and compliance with contact damping $\tau_{\tc, i}$ and compliance $\epsilon_{\tc, i}$.
Further, $\vec{\lambda}_\tc$ are the Lagrange multipliers corresponding to impulses normal to the contact surface, and $\vec{w}_\tc$ are so-called \emph{slack variables} that only play a role in defining the problem.
It is worth noting that since collisions typically happen at time-scales much faster than the time-step $h$, solving this system is not enough to handle collisions properly.
A separate impulse propagation step is needed, but this is beyond the current scope.
Also note that \eqref{eq:spook_inequality} is no longer a linear system but rather falls under the category of so-called \emph{mixed linear complementarity problems} (MLCPs).
Solving $\eqref{eq:spook_inequality}$ requires specialized methods \cite{renouf20053d}.

It is also possible to formulate Coulomb friction as an MLCP.
The Coulomb friction force $\vec{\mathrm{\vec{f}}}_{\text{Coulomb}}$ is proportional to the normal force and acts in the direction to oppose the sliding of two surfaces, either as a result of rotation or translation of the associated bodies:
\begin{align}
\label{eq:coulomb_linear}
    \| \vec{\mathrm{f}}_{\text{Coulomb}}\| \leq \mu\|\vec{\mathrm{f}}_n\|.
\end{align}
For rotational degrees of freedom, the corresponding inequality is
\begin{align}
    \label{eq:coulomb_rotational}
    \| \vec{\tau}_{\text{Coulomb}}\| \leq r \mu_r\|\vec{\mathrm{f}}_n\|,
\end{align}
where $r$ is the distance between the contact point and the body's center of mass.
In the MLCP model, the complementarity conditions enforce box constraints for the tangential constraint forces coupled to the normal multipliers, meaning that we use linearized norms as compared to \eqref{eq:coulomb_linear} and \eqref{eq:coulomb_rotational}.
We consider contact constraints $\vec{c} = \begin{bmatrix} \vec{c}_1^T \hdots \vec{c}_{n_\tc}^T \end{bmatrix}^T$, where $\vec{c}_i: \mathcal{M} \to \mathbb{R}^{d_{\tc, i}}$ has the Jacobian matrix $\mat{C}_i$.
Now, take $\mat{\bar{C}}_{i, j}: \mathcal{M} \to \mathbb{R}^{1 \times 6n_b}$ for $j \in \{1,\hdots,n_{\text{proj}, i}\}$ to be matrices that project vectors to the tangent of contact constraint $i$, and such that there is always a $j$ for which $\mat{\bar{C}}_{i, j}\vec{v} \geq \vec{0}$.
Finally, put $\bar{\mat{D}}_i = \begin{bmatrix} \bar{\mat{C}}_{i, 1}^T & \hdots & \bar{\mat{C}}_{i, n_{\text{proj}, i}}^T\end{bmatrix}^T$. 
With multiple frictional contact constraints, the equation for the velocity update reads
\begin{equation}
    \label{eq:spook_dry}
    \begin{aligned}
        \begin{bmatrix}
            \mat{M}_k       & -\mat{C}_k^T   & -\bar{\mat{D}}_k^T & \mat{0}      \\
            \mat{C}_k       & \mat{\Sigma}_\tc & \mat{0}            & \mat{0}      \\
            \bar{\mat{D}}_k & \mat{0}        & \mat{\Sigma}_\text{t}     & \mat{E}      \\
            \mat{0}         & \mat{U}        & -\mat{E}^T         & \mat{\Delta}
        \end{bmatrix}
        \begin{bmatrix}
            \vec{v}_{k+1} \\
            \vec{\lambda}_\tc      \\
            \lambdat      \\
            \vec{\sigma}
        \end{bmatrix}
         & +
        \begin{bmatrix}
            -\mat{M}\vec{v}_k + h (\vec{f} + \vec{f}_g)                         \\
            \frac{4}{h}\mat{\Upsilon} \vec{c}_k - \mat{\Upsilon} \vec{c}_k \vec{v}_k \\
            \vec{0}                                                                   \\
            \vec{0}
        \end{bmatrix}
        =
        \begin{bmatrix}
            \vec{0}   \\
            \vec{w}_\text{n} \\
            \vec{w}_\text{t} \\
            \vec{w}_\sigma
        \end{bmatrix},                                 \\
        \complementarityalign{\lambdan}{\vec{w}_\text{n}}, \\
        \complementarityalign{\lambdat}{\vec{w}_\text{t}}, \\
        \complementarityalign{\vec{\sigma}}{\vec{w}_\sigma},
    \end{aligned}
\end{equation}
where $\mat{U}$ is a block-diagonal matrix of friction coefficients, $\mat{\Delta}$ is a regularization matrix, $\mat{E}$ is a block-diagonal matrix of ones, and $\vec{\sigma}$ are multipliers that act to limit the tangential force.
Note that the block-antisymmetry is broken by $\mat{U}$ not having an antisymmetric partner. The block-diagonal entries of $\mat{U}$ are $\mat{U}_1, \hdots, \mat{U}_{n_\tc}$. Each $\mat{U}_i$ is a row vector or scalar with friction coefficients $\mu$ or $r \mu_r$.
For a detailed derivation and motivation of \eqref{eq:spook_dry}, we refer the reader to \cite{lacoursiere2007ghosts}.
We write \eqref{eq:spook_dry} compactly as
\begin{equation}
    \label{eq:spook_dry_compact}
    \begin{aligned}
        \begin{bmatrix}
            \mat{M}_k & -\jacb_k^T     \\
            \jacb_k & \mat{\Sigma}_\tB \\
        \end{bmatrix}
        \begin{bmatrix}
            \vec{v}_{k+1}   \\
            \vec{\lambda}_\tB \\
        \end{bmatrix}
         & +
        \begin{bmatrix}
            \vec{p} \\
            \vec{b} \\
        \end{bmatrix}
        =
        \begin{bmatrix}
            \vec{0}   \\
            \vec{w}_\tB \\
        \end{bmatrix}, \\
        \complementarityalign{\vec{\lambda}_\tB}{\vec{w}_\tB},
    \end{aligned}
\end{equation}
where $\mat{\Sigma}_\tB$, $\mat{B}_k$, $\vec{\lambda}_\tB$,  $\vec{b}$ and $\vec{w}_\tB$ are given through direct comparison.
Note that $\mat{\Sigma}_\tB$ is neither symmetric nor block-antisymmetric.
\subsection{Solver}
For a general multibody system with multiple bodies, multiple constraints, and dry friction, we need to solve for the velocity update $\vec{v}_{k+1}$ of the combined system \eqref{eq:spook_wihtout_complementarity} and \eqref{eq:spook_dry_compact}:
\begin{equation}
    \label{eq:mlcp}
    \begin{aligned}
        \begin{bmatrix}
            \mat{M}_k & -\jaca_k^T     & -\jacb_k^T   \\
            \jaca_k & \mat{\Sigma}_\tA & \mat{0}      \\
            \jacb_k & \mat{0}        & \drydiagonal \\
        \end{bmatrix}
        \begin{bmatrix}
            \vec{v}_{k+1}   \\
            \vec{\lambda}_\tA \\
            \vec{\lambda}_B \\
        \end{bmatrix}
         & +
        \begin{bmatrix}
            \vec{p} \\
            \vec{a} \\
            \vec{b} \\
        \end{bmatrix}
        =
        \begin{bmatrix}
            \vec{0}   \\
            \vec{0}   \\
            \vec{w}_b \\
        \end{bmatrix}, \\
        \complementarityalign{\lambdab}{\vec{w}_{\tB}}.
    \end{aligned}
\end{equation}
To do so, we start with a Schur complement operation of \eqref{eq:mlcp}, thus dividing the problem into two steps:
First, we solve for the multipliers $\lambdaa$, $\lambdab$ by
\begin{equation}
    \label{eq:mlcp_schur}
    \begin{aligned}
        \begin{bmatrix}
            \mat{A}_k \mat{M}^{-1}_k \mat{A}^T + \mat{\Sigma}_\tA & \mat{A}\mat{M}^{-1}\mat{B}^T                \\
            \mat{B}\mat{M}^{-1}\mat{A}^T                        & \mat{B}\mat{M}^{-1}\mat{B}^T + \drydiagonal \\
        \end{bmatrix}
        \begin{bmatrix}
            \lambdaa \\
            \lambdab \\
        \end{bmatrix}
         & +
        \begin{bmatrix}
            \vec{a} - \mat{A}\mat{M}^{-1} \vec{p} \\
            \vec{b} - \mat{B}\mat{M}^{-1} \vec{p} \\
        \end{bmatrix}
        =
        \begin{bmatrix}
            \vec{0} \\
            \vec{w}_{\tB}     \\
        \end{bmatrix}, \\
        \complementarityalign{\lambdab}{\vec{w}_{\tB}},
    \end{aligned}
\end{equation}
and then we solve for the velocity update by
\begin{equation}
    \label{eq:velocity_update}
    \mat{M}_k \vec{v}_{k+1} = \mat{M}\vec{v}_k + \mat{A}_k^T \lambdaa + \mat{B}_k^T \lambdab + h (\vec{f} + \vec{f}_g).
\end{equation}
In short, we solve \eqref{eq:mlcp_schur} by reducing it to a \emph{linear complementarity problem} (LCP) using block-LDLT factorization.
The LCP is solved using Lemke's algorithm \cite{cottle2009linear}, and forward and backward substitution completes the solution.
Details are presented in \ref{appendix:solver}.
This is the core of the custom physics simulator that is used in this work.

\subsection{Inverse Dynamics}
\label{sec:multibody_in_optimization}
In fields such as optimal control, reinforcement learning, and system identification, simulators are sometimes used in optimization loops  \cite{newburyReviewDifferentiableSimulators2024}.
Incorporating the dynamics (i.e., the time-stepper) in the optimization problem can be done in two ways:
First, through \emph{forward dynamics}---using the time-stepper to predict the next state or sequence of states. Second, by reformulating the time-stepper for \emph{inverse dynamics}---finding the forces required to realize a given state transition or trajectory.
For \eqref{eq:mlcp}, the external force $\vec{f}_k$ required to go from velocity $\vec{v}_k$ to velocity $\vec{v}_{k+1}$ is
\begin{equation}
    \vec{f}_k = h^{-1} \mat{M} (\vec{v}_{k+1} - \vec{v}_k) - h^{-1} \mat{A}^T \vec{\lambda}_\tA - h^{-1} \mat{B}^T \vec{\lambda}_\tB - \vec{f}_{g,k},
\end{equation}
where $\vec{\lambda}_\tA$ is given by
\begin{equation}
    \mat{\Sigma}_\tA \vec{\lambda}_\tA = \vec{b} - \vec{A}_k \vec{v}_{k+1}, \\
\end{equation}
and $\vec{\lambda}_\tB$ are given by
\begin{equation}
    \begin{aligned}
        \mat{\Sigma}_\tB \vec{\lambda}_\tB &= - \mat{B}_k \vec{v}_{k+1} + \vec{w}_\tB, \\
        \complementarityalign{\vec{\lambda}_\tA}{\vec{w}_\tB}.
    \end{aligned}
\end{equation}

\subsection{Automatic differentiation}
To use gradient-based optimization methods, it is necessary to take derivatives through the forward or inverse dynamics.
Most phyics engines do not provide such funcionality and one has to resort to finite differences or other black-box methods.
While exposing analytical derivatives in the API is an option, this method lacks flexibility.
Another way to achieve this is to write the physics simulator using a differentiable programming paradigm, usually with a technique known as automatic differentiation (AD).

AD is a family of methods for computing derivatives using the chain rule of calculus \cite{baydinAutomaticDifferentiationIn2018}.
Code written in a language or library that supports AD allows derivative information to be passed alongside the standard variables in function calls.
This is achieved under the hood by tracking every operation, knowing the derivatives of elementary operations, and applying the chain rule.

Different versions (or \emph{modes}) of AD appear as there is freedom in which order to compute the products in the chain rule.
We typically talk about two choices, although a mix also possible.
The first is \emph{forward mode accumulation} where so-called \emph{tangents} are passed alongside the standard variables or \textit{primals}. In forward mode accumulation, one typically considers derivatives with respect to a single scalar $x$ at a time.
Each elementary operation $\vec{\varphi}: \mathbb{R}^n \rightarrow \mathbb{R}^m$ is passed both an input primal $\vec{v}_{\text{in}} \in \mathbb{R}^n$ and an input tangent $\dot{\vec{v}}_{\text{in}} = \frac{d \vec{v}_{\text{in}}}{d x} \in \mathbb{R}^n$.
It is expected to return an output primal $\vec{v}_{\text{out}} \in \mathbb{R}^m$ and output tangent $\dot{\vec{v}}_\text{out} \in \mathbb{R}^m$:
\begin{align}
    \vec{v}_{\text{out}} &= \vec{\varphi}(\vec{v}_{\text{in}}), \\
    \dot{\vec{v}}_{\text{out}} &= \frac{\partial \vec{\varphi}}{\partial \vec{v}_{\text{in}}} \dot{\vec{v}}_{\text{in}}.
\end{align}
The tangent is given by a Jacobian-vector-product (JVP), which should have a known explicit formula.

The second mode is known as \emph{reverse mode accumulation} and requires a backward pass after all standard variables are known.
In the backward pass, so-called \emph{adjoints} are computed and passed in the reverse direction of the computational graph.

It is often stated that if a composite function $\mathbb{R}^n \to \mathbb{R}^m$ has many more inputs than outputs, $n \gg m$, then reverse mode accumulation is the most efficient. 
On the other hand, if $m \gg n$, forward mode accumulation is most efficient.
This is true if all operations are dense.
However, if the intermediate Jacobian matrices are sparse, there are typically ways to make use of sparsity.

While automatic differentiation is a powerful technique, it also has its limitations.
One can encounter exploding and vanishing derivatives due to finite numerical precision, resulting in unbounded errors.
AD has no mechanism for simplifying expressions (Cases where $\infty \cdot 0$ should be finite).
These cases must be identified and replaced with custom differentiation rules \cite{metzGradientsAreNot2022}.

The physics simulator used in this work is written in Python with the library \textsc{Jax} \cite{jax2018github}.
\textsc{Jax} provides a functional style of differentiating Python functions and does so by function tracing. The JVPs of all elementary operations are included in the library.
However, when implementing specific algorithms, it is more efficient to define custom rules.
We present the custom rules used in this work in \ref{appendix:jvp_rules}.

\section{State and parameter estimation}
\label{sec:parameter_identification}
In this section, we formulate and motivate the joint parameter and state estimation problem as a nonlinear least squares optimization problem.
We then motivate using inverse dynamics, the Levenberg-Marquardt optimization algorithm, and means to deal with rotations and exploit sparsity.
Finally, we address how we do linearized sensitivity analysis within this framework.

\subsection{The least squares optimization problem}
Let $\vec{x}_k = (\vec{q}_k, \vec{v}_k) \in \mathcal{M} \times (\mathbb{R}^6)^{n_\text{b}} = T\mathcal{M}$ denote the state of the multibody system at time-step $k$. 
Denote the parameters that we want to estimate by $\param \in B$ (mass, inertia, frictional parameters, etc.) where $B$ is a convex subset of $\mathbb{R}^{n_\text{p}}$ that excludes non-physical parameter values.
We will assume that the states evolve according to
\begin{equation}
    \label{eq:state}
    \vec{x}_{k+1} = \vec{\step}(\vec{x}_k, \vec{u}_k; \param) \boxplus \vec{\eta}_k,
\end{equation}
where $\vec{\step}: T\mathcal{M} \times \mathbb{R}^{n_\tc} \times B \to T\mathcal{M}$ is the time-stepper seen as a state transition function.
The term $\vec{\statenoise}_k$ is an additive \emph{state noise} and represents model errors and external disturbances affecting future states.
The vector $\vec{u}_k$ denotes the $n_\tc$-dimensional control signal at time $k$.
The $\boxplus$ operator was introduced in \eqref{eq:boxdef}, and it collapses to normal addition for velocity components.

We assume that the available data comes in the form of observations $\vec{y}_k \in \mathbb{R}^{n_\text{o}}$ that are made by an observation function $\vec{h}$:
\begin{equation}
    \label{eq:observation}
    \vec{y}_k  = \vec{h}(\vec{x}_{k}; \param^*) + \vec{\obsnoise}_k,
\end{equation}
where $\vec{\obsnoise}_k$ is an additive \emph{observation noise} and accounts for disturbances and errors in the observation function $\vec{h}$.

When identifying unknown parameters, we must choose how to deal with stochasticity.
Given some set of observation data, the most general solution is that of a multivariate probability distribution over the parameter space.
In this work, we reduce the complexity of the problem by assuming that the noise terms are zero-mean Gaussian, $\vec{\statenoise}_k \in \mathcal{N}(\vec{0}, \mat{Q}_k)$ and $\vec{\obsnoise}_k \in \mathcal{N}(\vec{0}, \mat{R}_k)$ for $k=1,\hdots,n$.
The consequence of this choice is that maximization of the likelihood function for both states and parameters can be formulated as a nonlinear least squares objective, see \ref{appendix:maximum_likelihood}.
The solvability of nonlinear least squares is the primary motivation behind this simple noise model and the chosen maximum likelihood approach.
We will return to how $\mat{Q}_k$ and $\mat{R}_k$ are specified.
For now, the maximum likelihood objective is that of finding the parameters and states that solve
\begin{mini}|l|
    {\param, \vec{x}_{0:n-1} \in B \times (T\mathcal{M})^{n}}{
    \sum_{i=0}^{n-1} \| \Delta \vec{x}_k \|^2_{\mat{Q}^{-1}_i} + \sum_{i=0}^{n-1} \| \Delta \vec{y}_k \|^2_{\mat{R}^{-1}_i}}
    {}{},
    \label{eq:mlo}
\end{mini}
where
\begin{equation}
    \vec{x}_{0:n-1} = (
        \vec{x}_0, \vec{x}_1, \hdots, \vec{x}_{n-1}
    ),
\end{equation}
the state residuals are given by
\begin{equation}
    \Delta \vec{x}_k(\param, \vec{x}_{0:n-1}) = \begin{cases}
        \vec{x}_{0}^{*} \boxminus \vec{x}_{0},  &k = 0, \\
        \vec{\step}(\vec{x}_{k-1}) \boxminus \vec{x}_{k}, \qquad &k = 1,\hdots, n-1,
    \end{cases}
\end{equation}
and the observation residuals are given by
\begin{equation}
    \Delta \vec{y}_k(\param, \vec{x}_k) = \vec{h}(\vec{x}_k; \param)  - \vec{y}_k, \qquad k = 0, \hdots, n - 1.
\end{equation}
We use the notation $\| \vec{x} \|_{\mat{P}}^2 = \vec{x}^T \mat{P} \vec{x}$, where $\mat{P}$ is positive definite, to denote weighted norms.
On top of the noise term assumptions, we also assume that we know the initial state to be $\vec{x_0}^*$ with some uncertainty as captured by a Gaussian distribution with zero mean and covariance matrix $\mat{Q}_0$.

\subsection{Weighing the residuals}
So far, we have made no assumption on the relation between the covariance matrices $\mat{R}_k$ and $\mat{Q}_k$ at different time-steps $k$.
Such assumptions are needed as there is no practical way to estimate all of $\mat{R}_k$ and $\mat{Q}_k$ separately.
For the observation noise, one can usually assume \emph{homoscedasticity}---that all $\vec{\obsnoise}_k$ have the same covariance $\mat{R}_k = \mat{R}_*$.
Moving forward, we follow this assumption and consider $\mat{R}_*$ to be a matrix of method parameters.

The same assumption should not be made for the state noise---we expect the uncertainty to depend on the local curvature of the kinematic constraints.
However, assuming that we make no numerical errors, the equations of motion only let errors enter as forces.
This is typically a good approximation even with some integration errors since the error in the force model is generally larger and can account for some integration errors \cite{sanjurjoAccuracyEfficiencyComparison2017}.
If we take this assumption seriously, we should replace the state errors with force or impulse errors.
Finding forces or impulses from a kinematic trajectory is precisely the inverse dynamics method described in section \ref{sec:multibody_in_optimization}.

That is, instead of computing the state residuals using forward dynamics, we opt to compute impulse residuals $\Delta \vec{p}_{k} = \begin{bmatrix}
    \Delta \vec{p}_{q,k}^T & \Delta \vec{p}_{u,k}^T
\end{bmatrix}^T$
using inverse dynamics:
\begin{equation}
    \Delta \vec{p}_{u, k+1}(\param, \vec{x}_{0:n-1}) = \mat{M} (\vec{v}_{k+1} - \vec{v}_k) - \mat{A}^T \vec{\lambda}_A - \mat{B}^T \vec{\lambda}_B - h \vec{\tau}_{g,k},
\end{equation}
\begin{equation}
    \Delta \vec{p}_{q, k+1}(\param, \vec{x}_{0:n-1}) = 
    \mat{M} (h^{-1}
    (\vec{\conf}_{k+1} \boxminus \vec{\conf}_{k}) 
    - \vec{v}_k) 
    - \mat{A}^T \vec{\lambda}_A -\mat{B}^T \vec{\lambda}_B - h \vec{\tau}_{g,k}.
\end{equation}
The residual $\Delta \vec{p}_{u, k+1}$ is the external impulse required to go from velocity $\vec{v}_k$ to velocity $\vec{v}_{k+1}$ and the residual $\Delta \vec{p}_{q, k+1}$ is the external impulse required to go from configuration $\vec{\conf}_k$ to configuration $\vec{\conf}_{k+1}$.

To make the connection to the notation we have developed so far, we should relate $\Delta \vec{p}_k$ to $\Delta \vec{x}_k$.
For the sake of compactness, we ignore complementarity conditions. The predicted velocity update $\vec{v}_{k+1}^* = \vec{\step}_u(\vec{x}_k)$ from state $\vec{x}_k$ is given by 
\begin{equation}
    \label{eq:no_impulse_noise_spook_update}
    \mat{M} \vec{v}_{k+1}^* = \mat{M} \vec{v}_k + \mat{A}^T \vec{\lambda}_A + h(\vec{\tau} + \vec{\tau}_g).    
\end{equation}
The corrected update that would get us to the next velocity $\vec{v}_{k+1}$ is
\begin{equation}
    \label{eq:impulse_noise_spook_update}
    \mat{M} \vec{v}_{k+1} = \mat{M} \vec{v}_k + \mat{A}^T \hat{\vec{\lambda}}_A + h(\Delta \vec{\tau}_{u,k+1}+ \vec{\tau} + \vec{\tau}_g),
\end{equation}
By subtracting \eqref{eq:no_impulse_noise_spook_update} from \eqref{eq:impulse_noise_spook_update}, we can formulate the impulse equation
\begin{equation}
    \begin{aligned}
        \begin{bmatrix}
            \mat{M}_k & -\jaca_k^T  \\
            \jaca_k & \mat{\Sigma}_A   \\
        \end{bmatrix}
        \begin{bmatrix}
            \Delta \vec{v}_{k+1}   \\
            \Delta \vec{\lambda}_A \\
        \end{bmatrix}
        =
        \begin{bmatrix}
            \Delta \vec{p}_{u,k} \\
            \vec{0} \\
        \end{bmatrix},
    \end{aligned}
\end{equation}
where 
\begin{equation}
    \Delta \vec{x}_k = \begin{bmatrix}
        \Delta \vec{\conf}_k \\
        \Delta \vec{v}_k
    \end{bmatrix}
    =
    \begin{bmatrix}
        \vec{\conf}_k^* \boxminus \vec{\conf}_k \\
        \vec{v}_k^* - \vec{v}_k
    \end{bmatrix}.
\end{equation}
Solving for $\Delta \vec{p}_{u,k+1}$ using Schur complements reveals the relation
\begin{equation}
    \label{eq:velocity_force_relation}
    \mat{M}_{\Sigma,k} \Delta \vec{v}_{k+1} = \Delta \vec{p}_{u,k+1},
\end{equation}
where $\mat{M}_{\Sigma,k} = \mat{A}_k^T \mat{\Sigma}^{-1}_A \mat{A}_k + \mat{M}_k$.

\begin{figure}
    \centering
    \includegraphics{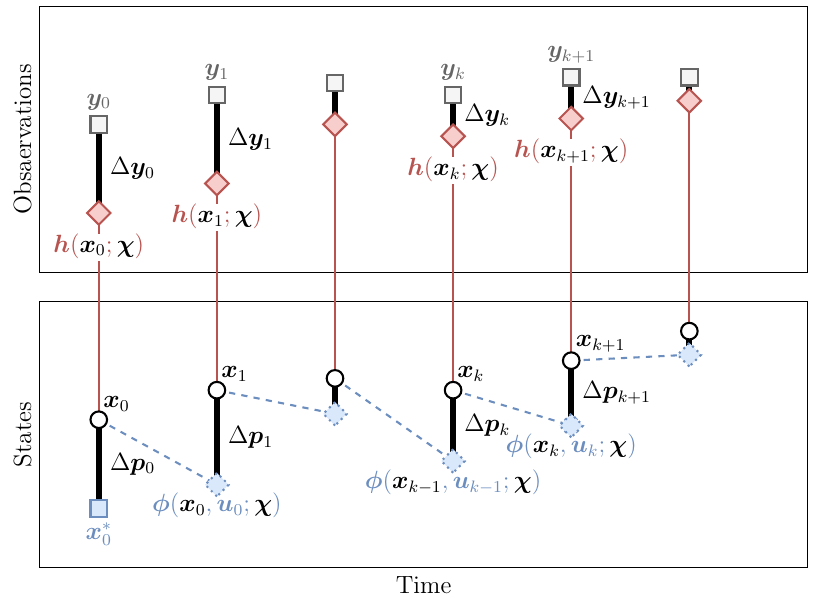}
    \caption{
        Abstract representation of the impulse residuals $\Delta \vec{p}_k$ and observation residuals $\Delta \vec{y}_k$.
        The white circles with black outlines represent the state optimization variables.
        These are used to predict the next states (blue diamonds) according to the discrete dynamics $\vec{\phi}$.
        However, the lines connecting them are dashed to indicate that this computation never takes place.
        Instead, the residuals $\Delta \vec{p}_k$ are computed directly using inverse dynamics.
        The state optimization variables also predict the observations (red diamonds).
        The observation residuals are given as the difference between the predicted observations and the measured observations (gray squares).
    }
    \label{fig:abstract_residuals}
\end{figure}

The same manipulations can be done for the configuration part of the state, revealing that
\begin{equation}
    \label{eq:conf_force_relation}
    h^{-1}\mat{M}_\Sigma \Delta \vec{\conf}_{k+1} =  \Delta \vec{p}_{q,k+1}.
\end{equation}
We are now able to state the optimization problem as
\begin{mini}|l|
    {\param, \vec{x}_{0:n-1} \in B \times (T\mathcal{M})^{n}}{
    \sum_{i=0}^{n-1} \| \Delta \vec{p}_{k} \|^2_{\mat{\Omega}_*^{-1}} + \sum_{i=0}^{n-1} \| \Delta \vec{y}_k \|^2_{\mat{R}^{-1}_*}}
    {}{},
    \label{eq:mlo2}
\end{mini}
where $\mat{\Omega}_*$ is the covariance matrix of the homoscedastic zero-mean impulse noise.
This is yet another matrix of method parameters.
An abstract representation of the residuals can be seen in \figref{fig:abstract_residuals}.

\subsection{State residuals at the initial time}
Finally, we should also form the residual corresponding to the state, $\Delta \vec{x}_0$.
It may not be immediately clear if this term is required or can safely be neglected.
However, we may expect that the constraint stabilization makes the outcome insensitive to this choice.
In that case, providing it has no downside, while not providing it could result in ill-conditioning.
On the other hand, if the outcome is sensitive to $\vec{x}_0$, it is only natural to provide a target state $\Delta \vec{x}_0^*$ and penalize deviations.
If no informed choice of $\vec{q}_0$ is available, we choose the target configuration $\vec{q}_0^*$ be the current guess $\vec{q}_0$ projected to the linearized holonomic constraint surface.
That is, the configuration residual is the distance to the constraint surface:
\begin{equation}
    \label{eq:residual_conf0}
    \Delta \vec{q}_0 = -(\mat{G}_0^T \mat{G}_0)^{-1} \mat{G}_0^T \vec{g}(\vec{q}_0).
\end{equation}
For the initial velocity, we let the target $\vec{v}_0^*$ be the current guess projected to the constraint surface tangent.
The residual is then the velocity normal to the constraint surface:
\begin{equation}
    \label{eq:residual_vel0}
    \Delta \vec{v}_0 = \mat{G}_0^T \mat{G}_0 \vec{v}_0.
\end{equation}
To get the impulse residuals, we multiply \eqref{eq:residual_conf0} and \eqref{eq:residual_vel0} by the projected mass matrix in accordance with \eqref{eq:conf_force_relation} and \eqref{eq:velocity_force_relation}.

\subsection{Solving nonlinear least squares}
The goal is to find the states $\vec{x}_{0:n-1}$ and parameters $\param$ that solve \eqref{eq:mlo}.
To make the notation more compact, we write the objective as
\begin{mini}|l|
    {\bar{\param} \in  B \times (T\mathcal{M})^{n}}{\| \combresiduals (\bar{\param}) \|^{2}_{\mat{W}}}
    {\label{eq:multiple_shooting_basic_compact}}{}.
\end{mini}
where $\bar{\param} = (\vec{x}_{0:n-1}, \param )$ combines the state trajectory and the free parameters, $\combresiduals$ is the concatenation of the impulse residuals and observation residuals,
\begin{equation}
    \combresiduals(\bar{\param}) 
    = \begin{bmatrix}
        \Delta \vec{p}_0(\bar{\param})^T & \hdots
        \Delta \vec{p}_{n-1}(\bar{\param})^T &
        \Delta \vec{y}_0(\bar{\param})^T & \hdots
        \Delta \vec{y}_{n-1}(\bar{\param})^T
    \end{bmatrix}^T,
\end{equation}
and $\mat{W}$ is the concatenation of the weight matrices
\begin{equation}
    \mat{W} = \text{diag}(\mat{\Omega}_0^{-1}, \hdots, \mat{\Omega}_{n-1}^{-1}, \mat{R}_0^{-1}, \hdots, \mat{R}_{n-1}^{-1}).
\end{equation}
Note that \eqref{eq:multiple_shooting_basic_compact} has a solution space that includes rotations and parameter limits.
This makes it a more general type of nonlinear least squares optimization problem.
To solve it, we use a Levenberg--Marquardt method \cite{levenbergMethodSolutionCertain1944}, for on-manifold optimization \cite{blanco2021tutorial} with box constraints.
In short, Levenberg--Marquardt is a trust-region extension of the Gauss--Newton method.
In Gauss--Newton, the parameters are iteratively updated by solving the local quadratic minimization problem. Levenberg--Marquardt adds regularization to ensure that every update progresses toward a local minima.
More precisely, we solve the following constrained linear least squares in each iteration:
\begin{mini}|l|
    {\vec{\delta}_{\bar{\param}}}{\norm{
        \combresiduals(\bar{\param}) + 
    \frac{\partial \combresiduals}{\partial \bar{\param}}
    \vec{\delta}_{\bar{\param}}
    }^2_{\mat{W}} + 
    \alpha \|\vec{\delta}_{\bar{\param}}\|^2,}
    {\label{eq:lm_update}}{}
    \addConstraint{\param + \vec{\delta}_{\param} \in B},
\end{mini}
where $\alpha$ is a damping factor chosen such that the cost is reduced
\begin{equation}
    \| \combresiduals(\bar{\param} \boxplus \vec{\delta}_{\param}) \|^2_{\mat{W}} < \| \combresiduals(\bar{\param}) \|^2_{\mat{W}}.
\end{equation} 
The parameters are updated by $\bar{\param} \leftarrow \bar{\param} \boxplus \vec{\delta}_{\bar{\param}}$. The pseudocode for the algorithm is shown in Algorithm \ref{alg:lm}.

\def\lmfunc{\Delta \vec{z}}
\def\lmvec{\overline{\lmfunc}}
\begin{algorithm}[ht!]
    \caption{Levenberg--Marquardt}\label{alg:lm}
    \algorithmicrequire $\lmfunc, \bar{\param}, \mat{W}, \varepsilon_g, \varepsilon_\delta$
    \begin{algorithmic}
        \State $\lmvec \gets\lmfunc(\bar{\param})$
        \State $\mat{J} \gets \texttt{jacfwd}(\lmfunc)(\bar{\param})$
        \State $\hat{\mat{H}}_{\lmfunc} \gets \mat{J}^T \mat{W} \mat{J}$
        \Comment{First order approximation of the Hessian}
        \State $\vec{g}_{\lmfunc} \gets \mat{J}^T \mat{W} \lmvec$
        \State $\alpha \gets \max(\text{diag}(\hat{\mat{H}})) \cdot 10^{-6}$
        \Comment{Damping term}
        \While {not \texttt{stop}}
        \State $\nu \gets 2$
        \State $\rho \gets 0$
        \While {$\rho \leq 0$ and $\alpha \leq 2^{32}$}
        \State $\vec{\delta}_{\bar{\param}} \gets -\text{SOL}(\hat{\mat{H}}_{\lmfunc} + \alpha \mat{I}, \vec{g}_{\lmfunc})$
        \State $\bar{\param}_p \gets \bar{\param} \boxplus \vec{\delta}_{\bar{\param}}$
        \Comment{Proposed update}
        \State $\rho \gets (\|\lmfunc(\bar{\param}_p)\|^2 - \norm{\lmvec}^2)/
            (\vec{\delta}_{\bar{\param}}^T (\alpha \vec{\delta}_{\bar{\param}} - \vec{g}_{\lmfunc}))$
        \Comment{Gain ratio}
        \If{$\rho \leq 0$}
        \State $\nu \gets \min(2\nu, 2^{32})$
        \State $\alpha \gets \min(\nu \alpha, 2^{32})$
        \EndIf
        \EndWhile
        \State $\texttt{stop} \gets \|\vec{g}_{\lmfunc}\|_2 \leq \varepsilon_g$ or $\|\vec{\delta}_{\bar{\param}}\|_2 \leq \varepsilon_\delta$ or $\rho \leq 0$
        \State $\bar{\param} \gets \bar{\param}_p$
        \State $\lmvec \gets\lmfunc(\bar{\param})$
        \State $\mat{J} \gets \texttt{jacfwd}(\lmfunc)(\bar{\param})$
        \State $\hat{\mat{H}}_{\lmfunc}\gets \mat{J}^T \mat{W} \mat{J}$
        \State $\vec{g}_{\lmfunc} \gets \mat{J}^T \mat{W} \vec{y}$
        \State $\alpha \gets \alpha \cdot \max(\frac{1}{3}, 1 - (2\rho - 1)^3)$
        \EndWhile
    \end{algorithmic} 
\end{algorithm} 
The method presented here shares similarities to specific Bayesian smoothing algorithms.
In particular, the so-called \emph{iterated extended Kalman smoother} is equivalent to the Gauss--Newton method to solve \eqref{eq:mlo} without tuning the parameters ($\param = \varnothing$), \cite{bellIteratedKalmanSmoother1994}.

\subsection{Computing rotation residuals}
The $\boxminus$ operator collapses to standard subtraction for all quantities except the rotational configuration.
To compute the generalized difference between two rotations, say $\vec{\quat}_a$ and $\vec{\quat}_b$, first form $\vec{\quat}' = \vec{\quat}^*_a \star \vec{\quat}_b$ with rotation vector representation $\vec{\psi}$. Let $e_s'$ be the scalar component of $\vec{e}'$ and $\vec{e}_v'$ be the vector component of $\vec{e}'$.
The explicit expression is
\begin{equation}
    \label{eq:exact_quat_residual}
    \vec{\quat}_a \boxminus \vec{\quat}_b = \vec{\rv} = \frac{2 \arccos (\quat_s') \vec{\quat}_v'}{\sqrt{1 - {\quat_s'}^2}}.
\end{equation}
Computing $\vec{\rv}$ from $\vec{\quat}'$ with this exact formula can be problematic as it has derivatives that approach infinity as $e_s'$ approaches $-1$.
We choose to avoid potential problems by replacing it with the small angle approximation
\begin{equation}
    \label{eq:rv_from_quat_approx}
    \vec{\quat}_a \boxminus \vec{\quat}_b = \vec{\rv} \approx 2 e_s' \vec{e}_v',
\end{equation}
when and only when we are computing rotation residuals. 
Note that other small angle approximations are possible since we expect to have $e_s' \approx 1$.

\subsection{Exploiting sparsity}
It is essential for a scalable method to utilize the sparsity pattern of the Jacobian $\frac{\partial \combresiduals}{\partial \bar{\param}}$ both during its construction and when solving the local linear least squares problem, \eqref{eq:lm_update}.
Now, if we write the Jacobian as a block matrix of sub-Jacobians
\begin{equation}
    \frac{\partial \combresiduals }{\partial \bar{\param}} = 
    \begin{bmatrix}
        \frac{\partial \Delta \vec{p}}{\partial \param} & \frac{\partial \Delta \vec{p}}{\partial \vec{x}_{0:n-1}} \\
        \frac{\partial  \Delta \vec{y} }{\partial \param} & \frac{\partial \Delta \vec{y}}{\partial \vec{x}_{0:n-1}}
    \end{bmatrix},
\end{equation}
it is easy to verify that the blocks $\partial \Delta \vec{p} / \partial \vec{x}_{0:n-1}$ and $\partial \Delta \vec{y} / \partial \vec{x}_{0:n-1}$ are block-diagonal and block-bidiagonal respectively.
This sparsity enables a drastic speed-up in the algorithmic differentiation, see \ref{appendix:jacobian} for our implementation, as well as in the box-constrained linear least squares solver, see \ref{appendix:box_constrainted_lls} for our implementation.

\subsection{Linearized parameter sensitivity}
Assuming that the observation noise is homoscedastic, $\mat{R}_* = \mat{R}_0 = \mat{R}_1 = \hdots = \mat{R}_{n-1}$, the unbiased sample covariance matrix of the observation noise is given by
\begin{equation}
    \hat{\mat{R}}_* = \frac{1}{n-1} \sum_{i=0}^{n-1} (\Delta \vec{y}_i) (\Delta \vec{y}_i)^T.
\end{equation}
The same is true for the impulse noise
\begin{equation}
    \hat{\mat{\Omega}}_* = \frac{1}{n-1} \sum_{i=0}^{n-1} (\Delta \vec{p}_{qu,i}) (\Delta \vec{p}_{qu,i})^T,
\end{equation}
from which it is easy to form $\hat{\mat{W}}$ is by
\begin{equation}
    \hat{\mat{W}} = \text{diag}(\hat{\mat{\Omega}}_0^{-1}, \hdots, \hat{\mat{\Omega}}_{n-1}^{-1}, \hat{\mat{R}}_0^{-1}, \hdots, \hat{\mat{R}}_{n-1}^{-1}).
\end{equation}
Linearizations of the combined state and observation residuals $\combresiduals$ around $\bar{\param}$ can be written as
\begin{equation}
    \combresiduals(\bar{\param} + \delta \bar{\param}) \approx 
    \combresiduals(\bar{\param}) +
    \frac{\partial \combresiduals}{\partial \bar{\param}} \delta \bar{\param}.
\end{equation}
The error approximately propagates as
\begin{equation}
    \hat{\mat{W}} = 
    \begin{bmatrix}
        \frac{\partial \combresiduals}{\partial \param} &
        \frac{\partial \combresiduals}{\partial \vec{x}_{0:n-1}}
    \end{bmatrix}
    \begin{bmatrix}
        \mat{\Lambda}_{\chi \chi} & \mat{\Lambda}_{\chi 0} \\
        \mat{\Lambda}_{0 \chi} & \mat{\Lambda}_{0 0}
    \end{bmatrix}
    \begin{bmatrix}
        \frac{\partial \Delta  \vec{z}}{\partial \param^T} \\
        \frac{\partial \Delta  \vec{z}}{\partial \vec{x}_{0:n-1}^T}
    \end{bmatrix},
\end{equation}
where $\vec{\Lambda}_{\chi \chi}$ is the covariance matrix of the parameters, $\vec{\Lambda}_{0 0}$ is the covariance matrix of the tunable states, and $\vec{\Lambda}_{\chi 0}$, $\vec{\Lambda}_{0 \chi}$ are cross-covariances.
It turns out that we can solve for $\mat{\Lambda}_{\chi \chi}$ by
\begin{equation}
    \label{eq:param_cov}
    \vec{\Lambda}_{\chi \chi} =  (\mat{T}^T \mat{T})^{-1} \mat{T}^T \hat{\mat{W}} \mat{T} (\mat{T}^T \mat{T})^{-1},
\end{equation}
where
\begin{equation}
    \mat{T} = \Big[
        \mat{I} - 
        \frac{\partial \combresiduals}{\partial \vec{x}_{0:n-1}}
        \Big(
            \frac{\partial \combresiduals}{\partial \vec{x}_{0:n-1}^T} \frac{\partial \combresiduals}{\partial \vec{x}_{0:n-1}}
        \Big)^{-1} 
        \frac{\partial \combresiduals}{\partial \vec{x}_{0:n-1}^T}
    \Big] \frac{\partial \combresiduals}{\partial \param},
\end{equation}
assuming that $\mat{T}$ and by extension $\partial \combresiduals / \partial \vec{x}_{0:n-1}^T$ has full column rank.
If not, the equations can be written in terms of QR-factorizations or held approximately with regularization.
Our implementation of the above uses regularization and a sparse linear solver to compute $\mat{T}$, followed by a singular value decomposition to compute $\mat{\Lambda}_{\chi \chi}$.
An interpretation of $\mat{T}$ is that it is the Jacobian corresponding to \emph{variable projection}, see \cite{o2013variable}, of the linearized state residuals.
The covariance matrix of the parameters contains all information for linearized sensitivity analysis.
A singular value decomposition of $\vec{\Lambda}_{\chi \chi}$ can be used for further analysis, such as determining parameter dependencies and null space.

\section{Experiments}
\label{sec:case_study}
The goal is to use the Levenberg-Marquardt algorithm with the techniques outlined in Section \ref{sec:parameter_identification} to identify unknown model parameters of a pendulum and a Furuta pendulum from recorded joint angles and control signal.
In this section, we investigate the performance and sensitivity to method parameters in both real and simulated scenarios.
The experiments are done using hardware from the PendCon Advanced series, the same hardware as used in \cite{leonid2017prop}, in both pendulum and Furuta pendulum configurations, see \figref{fig:furuta_images}.
Regularized time-discrete multibody models are the target for identification as described in Section \ref{sec:multibody_dynamics}.
The parameters subject to identification are moments of inertia, frictional parameters, direction of gravity, and motor gain.

\subsection{Hardware}
\begin{figure}
    \centering
    \includegraphics[height=8cm]{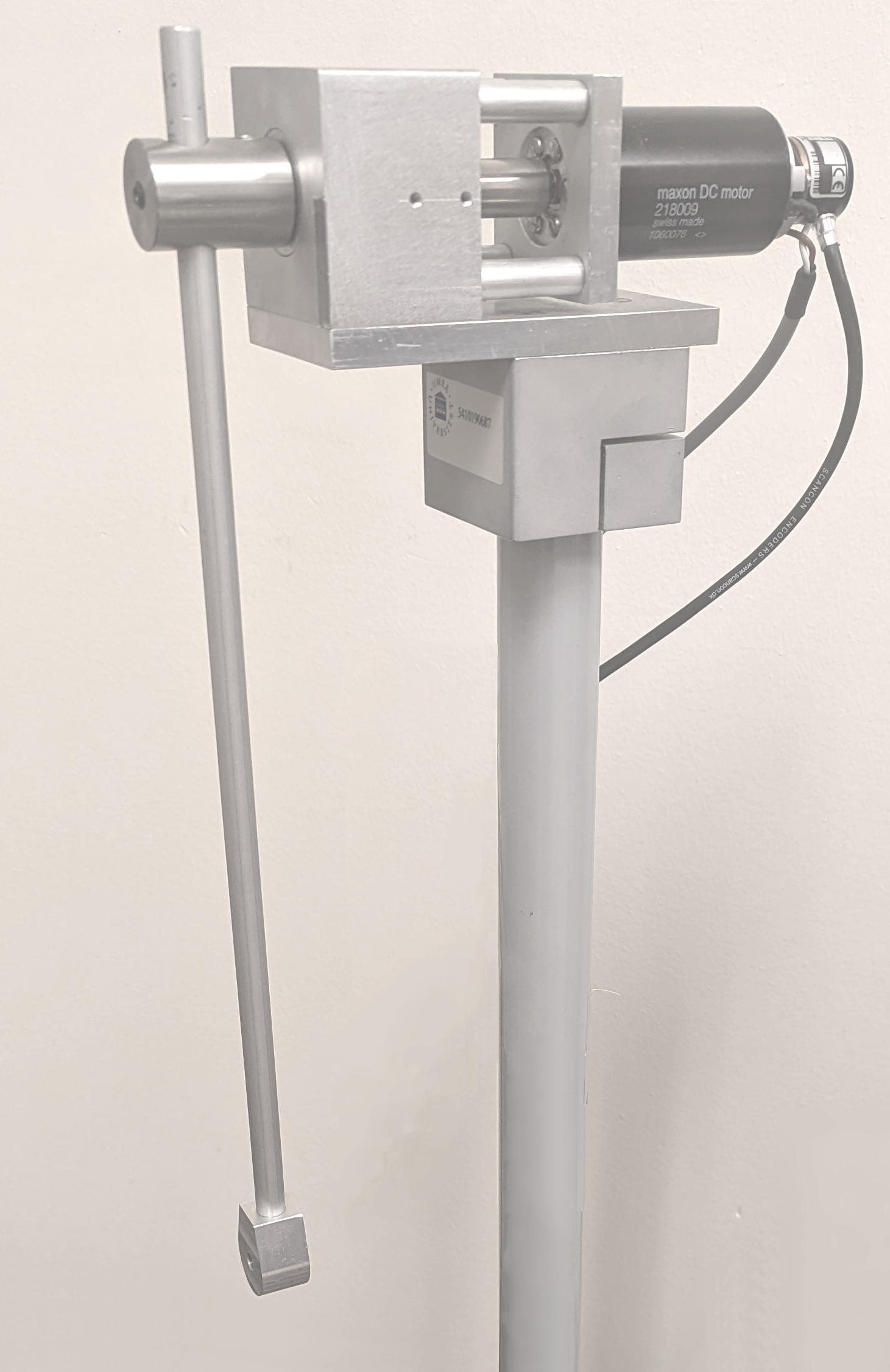}
    \includegraphics[height=8cm]{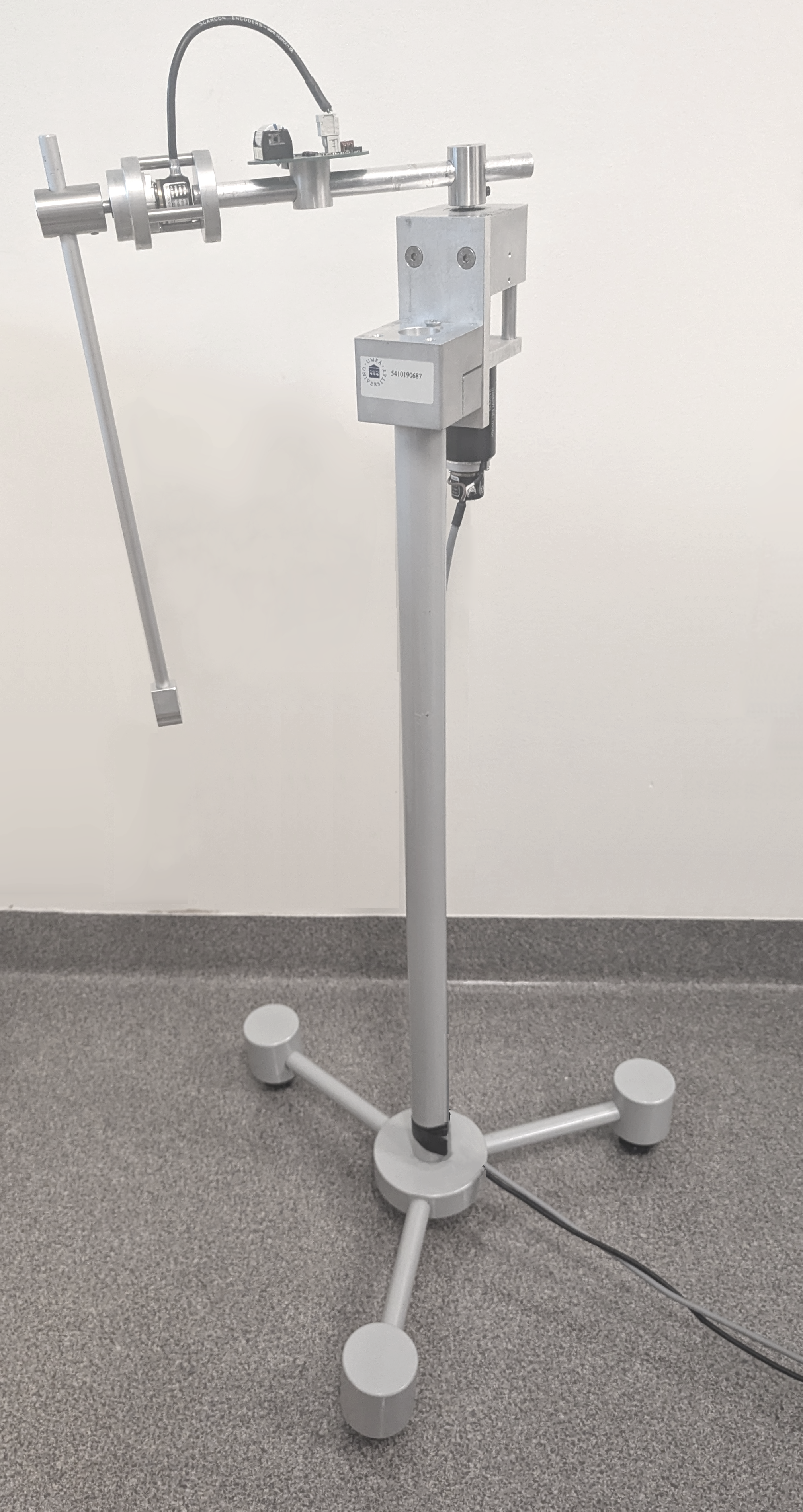}
    \caption{Images of the pendulum and the Furuta pendulum.}
    \label{fig:furuta_images}
\end{figure}

The PendCon Advanced series can be assembled as both a pendulum and a Furuta Pendulum, \figref{fig:furuta_images}.
The Furuta pendulum has two arms:
One arm is actuated by an electric motor and rotates about an axis fixed to the stand, while the other arm is passive and rotates about an axis fixed on the first arm.
Joint angles are measured by rotary encoders with a resolution of 4096 pulses per revolution.
The signal of the rotary encoder attached to the acuated joint is sent by wire, while the signal of the rotary encoder attached to the under-actuated joint is sent wirelessly.
To communicate with the hardware, we use \textsc{Matlab} Simulink Real-Time Explorer with a host PC and a dedicated target PC.
The target PC has a RS 232-card that talks to a conversion unit which receives the signals from the rotary encoders and sends a voltage (the control signal) to a PWM-amplifier that outputs a current proportional to that input voltage.
This current passes through the \textsc{Maxon} RE40 electric motor. We use a sample frequency of \SI{1000}{Hz} during all experiments.
Additional weight is put on the three feet to keep them fixed in place during experiments.

\subsection{Discrete pendulum}
The discrete pendulum consists of a single rigid body attached to the global frame with a \emph{frictional hinge joint}, a combination of kinematic constraints that restricts motion to rotation in a plane while also modeling friction, see \ref{appendix:hinge}.
The hinge is displaced a distance $l$ from the body's center of mass and set to have compliance $\SI{1e-4}{m/N}$ and damping \SI{0.02}{s}.
We label the mass $m_\tP$ and the inertia in the plane of oscillation $J_\tP$.
The hinge joint has a viscous friction coefficient $b$ and dry friction coefficient $\mu$. However, the dry frictional coefficient is coupled to the internal geometry such that we can only identify $r \mu$ where $r$ has the unit of distance.
The external force is given by gravity, $\vec{\mathrm{f}}_\tP = \begin{bmatrix} 0 & -m_\tP g_{acc} & 0 \end{bmatrix}^T$, where $g_{acc} = 9.82\,\si{m/s^2}$.
We label the angle with respect to the downright position in the plane of oscillation $\theta$.
We also consider a \emph{Stribeck friction} model that models the transitions between stick and slip.
If this model is used, the dry friction coefficient at time-step $k+1$ is given by
\begin{equation}
    \label{eq:stribeck}
    \mu_{k+1} = \mu_{\text{kinematic}} + \mu_{\text{delta}} \exp (-\zeta \omega_k^2),
\end{equation}
where $\mu_{\text{kinematic}}$ is a kinematic friction coefficient, $\mu_\text{delta} = \mu_\text{static} - \mu_\text{kinematic}$ is the (positive) difference between static and kinematic friction coefficients, $\zeta$ models the transition rate, and $\omega_k$ is the angular velocity in the plane of oscillation at time-step $k$.

\subsection{Discrete Furuta pendulum}
\begin{figure}
    \centering
    \includegraphics{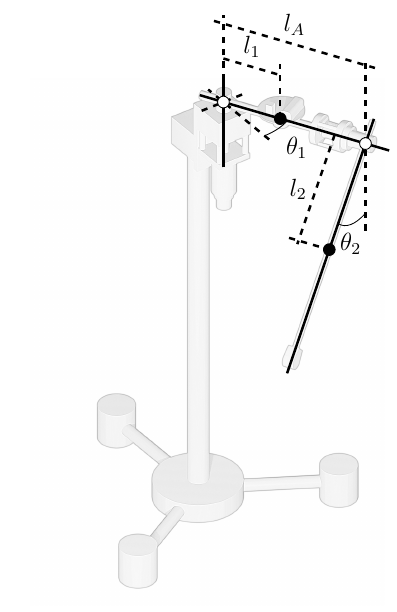}
    \raisebox{2.3\height}{
    \begin{minipage}[t]{.45\textwidth}%
    \includegraphics{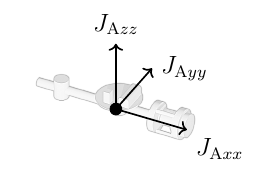}
    \includegraphics{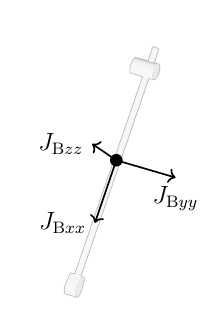}
    \end{minipage}}
    \caption{
        A render of the Furuta pendulum with annotated joint angles and lengths. 
        The two rigid bodies are drawn separately to the right, and their respective principal components of inertia are indicated.
    }
    \label{fig:furuta_labeled}
\end{figure}

The simulation model of the Furuta pendulum consists of two rigid bodies and two frictional hinge joints.
The stand is assumed not to flex or move and is thus not modeled.
We label the body that connects to the global frame $\tA$ and the other body $\tB$.
We index the hinge that connects $\tA$ to the global coordinate system by $1$ and the hinge that connects $\tA$ and $\tB$ by $2$.
The parameters $l_1$, $l_2$, and $l_A$ specify the offset between the center of mass for $\tA$ and $\tB$ and the hinges as indicated in \figref{fig:furuta_labeled}.
Hinge 1 has the viscous friction coefficient $b_1$ and the dry friction coefficient $\mu_1$ associated with the internal distance $r_1$.
Similarly, hinge 2 has the viscous friction coefficient $b_2$, the dry friction coefficient $\mu_2$, and the internal distance $r_2$.
We distinguish between the linear compliances $\epsilon_1^{(\ell)}$ and $\epsilon_2^{(\ell)}$, and the rotational compliances $\epsilon_1^{(r)}$ and $\epsilon_2^{(r)}$.
The external force is given by gravity as $\vec{\mathrm{f}}_\tA = m_A \vec{g}_{\text{acc}}$, $\vec{\mathrm{f}}_\tB= m_B \vec{g}_{\text{acc}}$ where $\norm{\vec{g}_{\text{acc}}} = 9.82\,\si{m/s^2}$ is given but the direction is a model parameter.
We also model the fictitious forces by $\vec{\tau}_{g, \tA} =\angvel_\tA \times \mat{J}_\tA(\vec{\conf}_\tA) \angvel_\tA$ and similarly for body $\tB$.

An ideal DC motor develops a torque proportional to the current passing through the stator coil, while the shaft's angular velocity induces a voltage that opposes the change in current that induced it.
In the PendCon Advanced series, the PWM amplifier adjusts the voltage (up to \SI{20}{V}) such that the current passed through the DC motor is proportional to the input voltage $u_{\text{in}}$.
We assume no delay or artifacts from the digital to analog converter, which allows us to directly associate $u_{\text{in}}$ with the control signal, $u = u_{\text{in}}$.
Thus, we model the torque as proportional to the control signal $\tau_{\text{motor}} = K u$ where $K$ is the motor gain.
To account for the extra inertia of the shaft, $J_{\text{shaft}}$, we model the motor through a regularized nonholonomic constraint, see \ref{appendix:hinge}.

The masses of $\tA$ and $\tB$ are $m_\tA$, $m_\tB$ and their moments of inertia are assumed to have principal axes aligned to their local frames:
\begin{eqnarray}
    \mat{J}_A^{(\tA)} = \text{diag}(J_{\tA xx}, J_{\tA yy}, J_{\tA zz}), \nonumber \\
    \mat{J}_\tB^{(\tB)} = \text{diag}(J_{\tB xx}, J_{\tB yy}, J_{\tB zz}).
\end{eqnarray}
All fixed parameter values and initial parameter values subject to optimization are shown in Table \ref{tab:furuta_default_parameters}.
Note that reason we put $r_1 = r_2 = 1\,\si{m}$ is not physical.
Instead, this is an arbitrary choice we make to unify friction due to force and torque multipliers.

\begin{table}
\caption{
    Default parameters of the discrete Furuta pendulum. 
    The left table contains fixed parameter values, and the right table contains initial parameter values that are subject to identification.
}
\label{tab:furuta_default_parameters}
\begin{tabular}{l|rrrrr|r}
    \midrule
    $l_1$ & \SI{0.128}{m} \\
    $l_\tA$ & \SI{0.248}{m} \\
    $l_2$ & \SI{0.92}{m} \\
    $m_\tA$ & \SI{0.238}{kg} \\
    $m_\tB$ & \SI{0.428}{kg} \\
    $J_{\tB xx}$ & $10^{-4}\,\si{kg m^2}$ \\
    $J_{\text{shaft}}$ & \SI{4e-4}{kg m^2} \\
    $r_1, r_2$ & \SI{1}{m} \\
    $\tau_1, \tau_2$ & $0.02\,\si{s}$ \\
    $\epsilon_1^{(\ell)}, \epsilon_2^{(\ell)}$ & $10^{-4}\,\si{m/N}$ \\
    $\epsilon_1^{(r)}, \epsilon_2^{(r)}$ & $10^{-4}\,\si{rad/N}$ \\
    \bottomrule
    \end{tabular}
    \qquad
    \begin{tabular}{l|rrrrr|r}
    \midrule
    $\diag{(\mat{J}_{A})}$ & $\begin{bmatrix}
        0.01 & 0.01 & 0.01
    \end{bmatrix}^T\,\si{kg m^2}$ \\
    $\diag{(\mat{J}_{B})}$ & $\begin{bmatrix}
        - & 0.01 & 0.01
    \end{bmatrix}^T\,\si{kg m^2}$ \\
    $\vec{g}_{\text{acc}}$ & $\begin{bmatrix}
        0 & -9.82 & 0
    \end{bmatrix}^T$\,\si{ms^{-2}} \\
    $b_{1}$ & $10^{-4}\,\si{kg s^{-1}}$ \\
    $b_{2}$ & $10^{-4}\,\si{kg s^{-1}}$ \\
    $\mu_{1}$ & $10^{-4}$ \\
    $\mu_{2}$ & $10^{-4}$ \\
    $K$ & \SI{0.1}{Nm V^{-1}} \\
    \bottomrule
    \end{tabular}
\end{table}

\subsection{Initialization}
\label{sec:initialization}
Before optimization, a set of initial states and parameters have to be chosen.
To select the initial states, we first smoothen the observation data using a low-pass filter.
The filter is a second-order Butterworth filter with \SI{10}{Hz} as the cutoff frequency.
It is applied twice to remove lag:
First in the forward direction and then in the reverse direction.
The filtered signal is downsampled to the desired timestep of the physics simulator.
The initial configurations $\vec{q}_{0:k-1}$ are chosen to satisfy the observation function at zero constraint violation.
The generalized velocities are approximated from the initial configurations using
\begin{equation}
    \vec{u}_{k+1} = \frac{1}{h}(\vec{q}_{k+1} \boxminus \vec{q}_{k}),
\end{equation}
where $h$ is the time-step.

To reduce the size of the hyperparameter space, we make the choice to do all parameter tuning with $\mat{R}_*^{-1} = \mat{I}$ and $\mat{\Omega}_*^{-1} = \kappa \mat{I}$, where $\kappa$ is a factor that we refer to as the \emph{state error weight}.

\subsection{Calibrating the pendulum}
We collect five time-series of observations where we release the pendulum from an almost upright position and let the oscillations die out, see \ref{appendix:pendulum}.
We consider optimization on three sets of parameters:
A \emph{viscous model} where $J$ and $b$ are optimized, a \emph{dry model} where $J_\tP$, $b$ and $r\mu$ are optimmized, and a \emph{Stribeck model} where \eqref{eq:stribeck} is used and all of $J_\tP$, $b$, and $r \mu_s$, $r \mu_{s-k}$ and $\zeta$ are optimized.
The model is tuned with respect to each time-series individually.
To compare the results, we compute a cross-validation metric.
More specifically, each tuned model is used to do state estimation of every time-series it was not tuned on.
The state estimation is done with the same method but with $\param = \varnothing$.
The cross-validation score is computed as the average cost of the state estimation.
All cost and cross-validation (CV) scores are reported as mean squared errors---the squared residuals divided by the number of time-steps.

The estimated parameters and the linearized standard deviation estimates of the dry pendulum model are shown in Table \ref{tab:param_dry}.
The inertia maximally varies with 0.8\%, a small relative error that indicates that the estimated parameters remain close to consistent between the time-series.
However, this consistency fails if the linearized standard deviation estimates are taken seriously---$10$ standard deviations are required to go from the inertia of time-series 0 to the inertia of time-series 3.
The average cost and cross-validation score are shown in Table \ref{tab:pendulum_error} for the three models.
The cross-validation score is close to the cost for all models (maximally varies by a factor $1.04$), likely explained by the similarity in excitation.
The cost is reduced by a factor $0.77$ from the viscous model to the dry model, suggesting that dry friction is important.
An additional reduction, almost as significant, is found from the dry model to the Stribeck model (a factor $0.80$).

\begin{table}
\caption{
    Estimated parameters of the pendulum.
    Values in parenthesis are the standard deviation estimates to the same precision as the least significant digit.
    If the estimated parameter value is $0$, then the value in parenthesis is the standard deviation in absolute precision.
}
\label{tab:param_dry}
\resizebox{\columnwidth}{!}{
\begin{tabular}{l|rrrrr|r}
\toprule
{} & {0} & {1} & {2} & {3} & {4} & {Combined} \\
\midrule
Cost & 2.06e-06 & 1.03e-06 & 2.04e-06 & 9.82e-07 & 1.25e-06 & 1.47e-06 \\
CV & 1.47e-06 & 1.49e-06 & 1.50e-06 & 1.47e-06 & 1.47e-06 & 1.48e-06 \\
$J\:(\si{\kilogram\meter^2})$ & 0.003916(2) & 0.003889(3) & 0.003912(3) & 0.003886(3) & 0.003888(4) &  \\
$b\:(\si{\kilogram\per\second})$ & 0(3e-05) & 0.00025(3) & 0.00026(4) & 0.00014(4) & 0.00020(4) &  \\
$r \mu\:(\si{m})$ & 0.00220(5) & 0.00173(5) & 0.00171(6) & 0.00189(6) & 0.00182(6) & \\
\bottomrule
\end{tabular}
}
\end{table}
\begin{table}
\caption{Comined errors of all pendulum models.}
\label{tab:pendulum_error}
\centering
\begin{tabular}{l|rrrrr|r}
\toprule
{} & {Viscous} & {Dry} & {Stribeck} \\
\midrule
Cost &  1.92e-06 & 1.47e-06 & 1.18e-06  \\
CV &  1.92e-06 & 1.48e-06 & 1.20e-06  \\
\bottomrule
\end{tabular}
\end{table}

\subsection{Calibrating the Furuta pendulum}
We perform three real-world Furuta pendulum experiments based on different excitation patterns (scenarios), resulting in three different trajectories that we use for parameter estimation.
The first scenario is the \emph{release scenario}, where the pendulum is released from its upright position.
The second scenario is the \emph{swing-up scenario}, where the pendulum is swung up from its downright position.
The third scenario is the \emph{pulse scenario}, where the electric motor is subject to a pulse control signal.
We present details about the data and how it is generated in \ref{appendix:data_real_furuta_pendulum}.

One gain, five inertial, four frictional parameters, and the stand's tilt are tuned.
We represent the tilt as a rotation applied to the gravitational acceleration vector.
As with the rotational part of the states, it is represented by a quaternion, resulting in a 13-dimensional parameter space.
The optimization is done with $\kappa = 100$ and $h = 0.01\,\si{s}$.
The iteration limit is set to $20$ and the tolerances are $\epsilon_g = 10^{1}$ and $\epsilon_\delta = 10^{-6}$.
As with the pendulum, we compute a cross-correlation metric by fitting each tuned model to the trajectories of the other two scenarios while only tuning the states.
Since the motor gain $K$ is never excited in the release scenario, we put it to $K = 7.0$ during cross-validation.

Results for parameter estimation of the Furuta pendulum are shown in Table \ref{tab:furuta_baseline_mix}.
The cross-validation score is the smallest for the model that is tuned on the swing-up scenario, \SI{2.34e-5}{}, closely followed by the release scenario, \SI{2.46e-5}{}.
The pulse scenario has an order of magnitude higher cross-validation score, \SI{2.36e-4}{}, indicating that it fails to excite the parameters to the same degree.
The parameter $J_{\tB yy}$ can be identified with the smallest relative error (as measured by parameter value over standard deviation estimate).
On the other hand, the parameters $J_{\tA xx}$ and $J_{\tA yy}$ have standard deviation estimates that are larger than the estimated values.
In practice, we lack information to specify $J_{\tA xx}$ and $J_{\tA yy}$.
This is explained by the fact that body $\tA$ never rotates about those axes.
Most parameters fail to remain consistent with the linearized parameter sensitivity over all three scenarios.
Still, the comparable cross-validation scores for swing-up and release suggest only minor differences during simulation.

\begin{figure}
    \centering
    \includegraphics[width=0.495\textwidth]{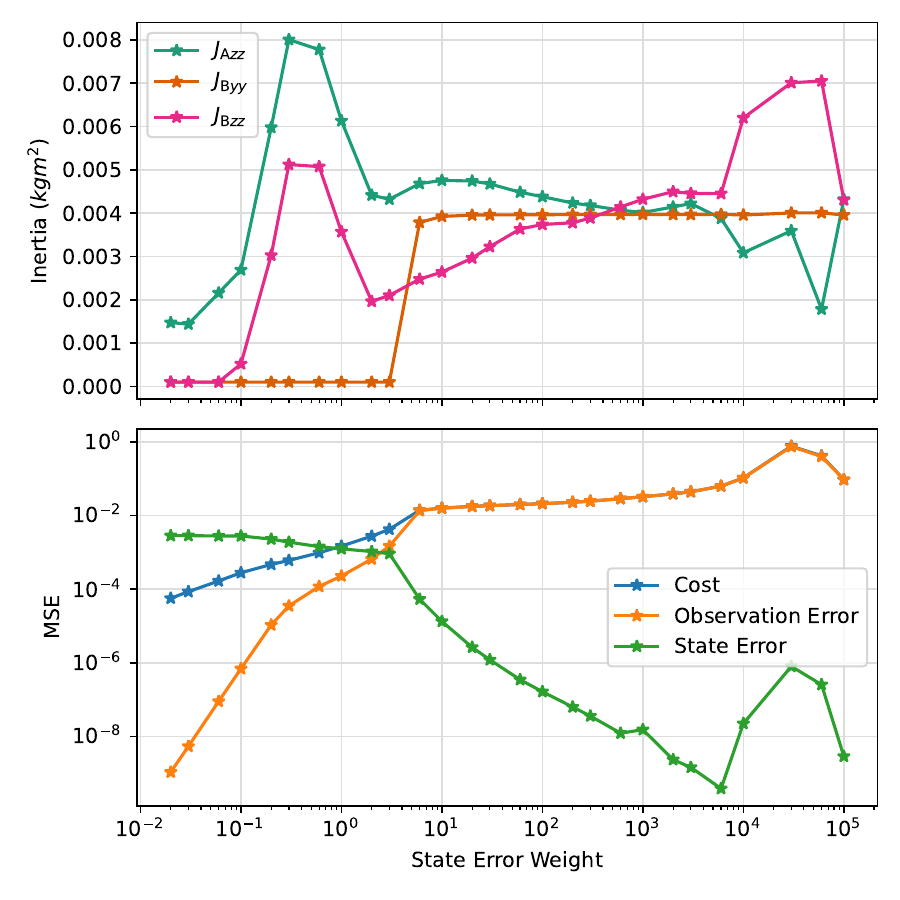}
    \includegraphics[width=0.495\textwidth]{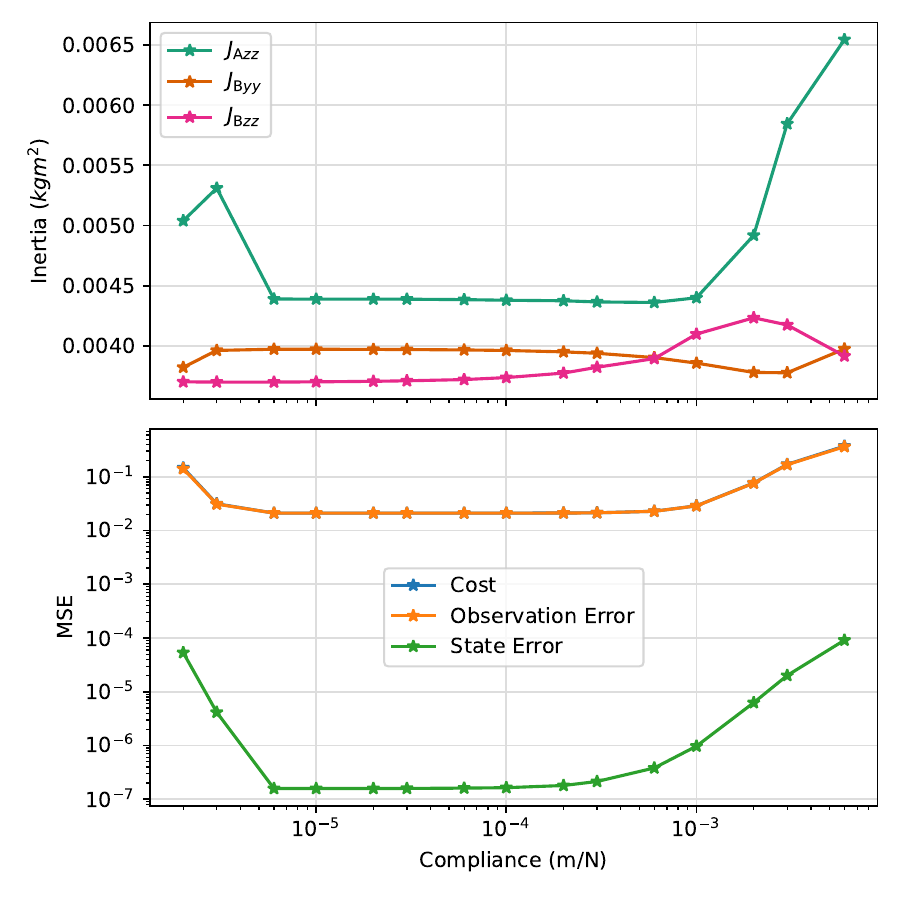}
    \caption{
        The two upper plots show three of the estimated parameters, tuned to the swing-up scenario versus different method parameters:
        State error weight in the left plot and compliance in the right plot.
        The lower plots show the mean squared error (MSE), plotted for observation error and state error versus the $\stateweight$ and $\epsilon$.
    }
    \label{fig:state_error_weight_sweep_real}
\end{figure}

Next, we focus on the swing-up scenario and the optimization result for different state error weights and global hinge compliance $\epsilon = \epsilon_1^{(\ell)} = \epsilon_2^{(\ell)} = (1\,\si{rad/N}) \epsilon_1^{(r)} = (1\,\si{rad/N}) \epsilon_1^{(r)}$. For the sake of comparison, we optimize over 20 iterations without other stopping criteria.
\figref{fig:state_error_weight_sweep_real} shows the estimated inertia components and cost.
The model overfits the real observations for $\stateweight \lesssim 10$, resulting in rapid variations of the identified parameters and high state error.
For $\kappa \gtrsim 5000$, the downward trend of the state error stops, and the estimated $J_{\tA zz}$ and $J_{\tB zz}$ become more sensitive to $\kappa$.
As for the hinge compliance, we observe an increase in both state error and observation error when the hinge is too soft, $\epsilon_1 \gtrsim 5 \times 10^{-4}\,\si{m/N}$, or too stiff, $\epsilon \lesssim 10^{-5}\,\si{m/N}$.

\begin{figure}
    \centering
    \includegraphics[width=0.495\textwidth]{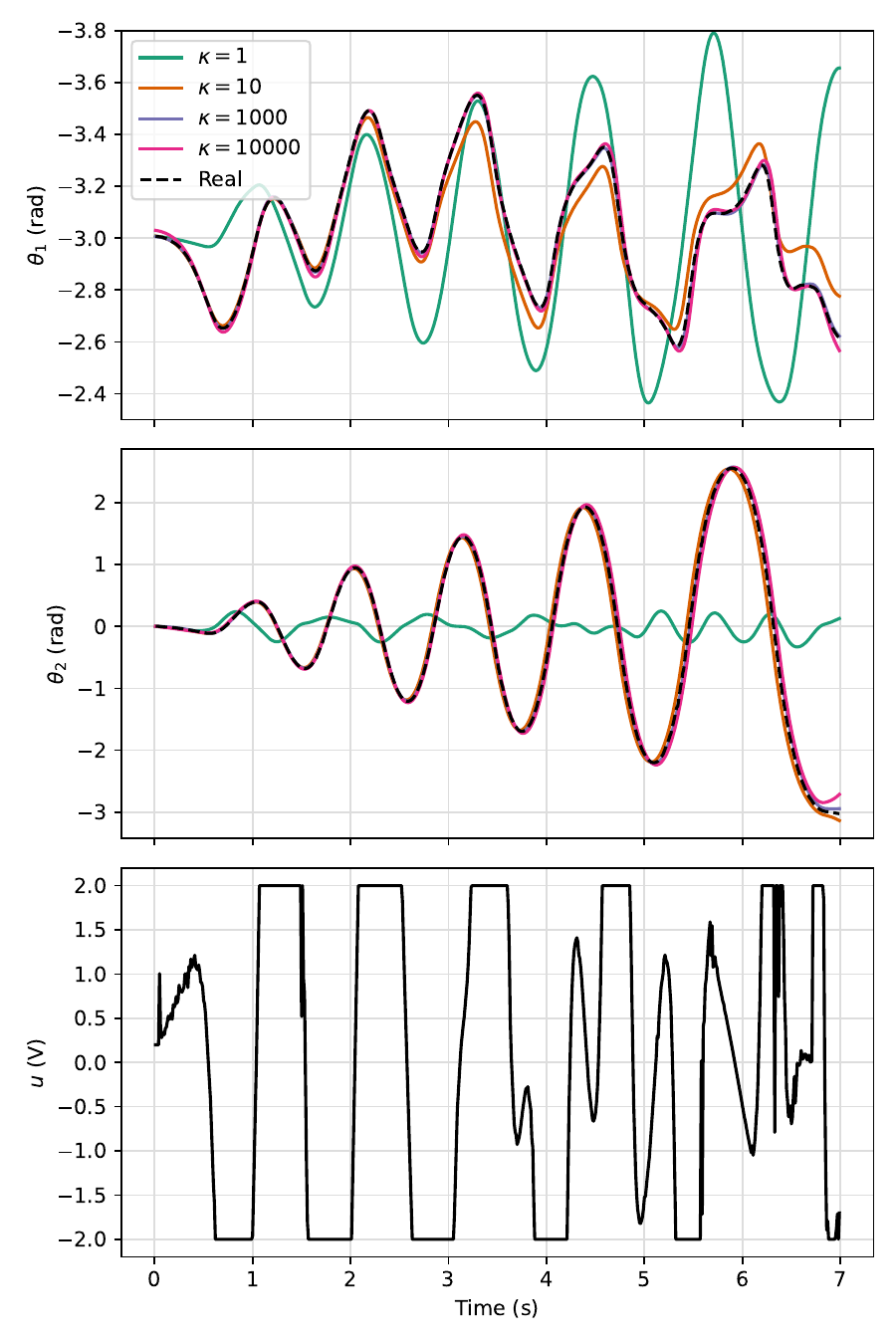}
    \includegraphics[width=0.495\textwidth]{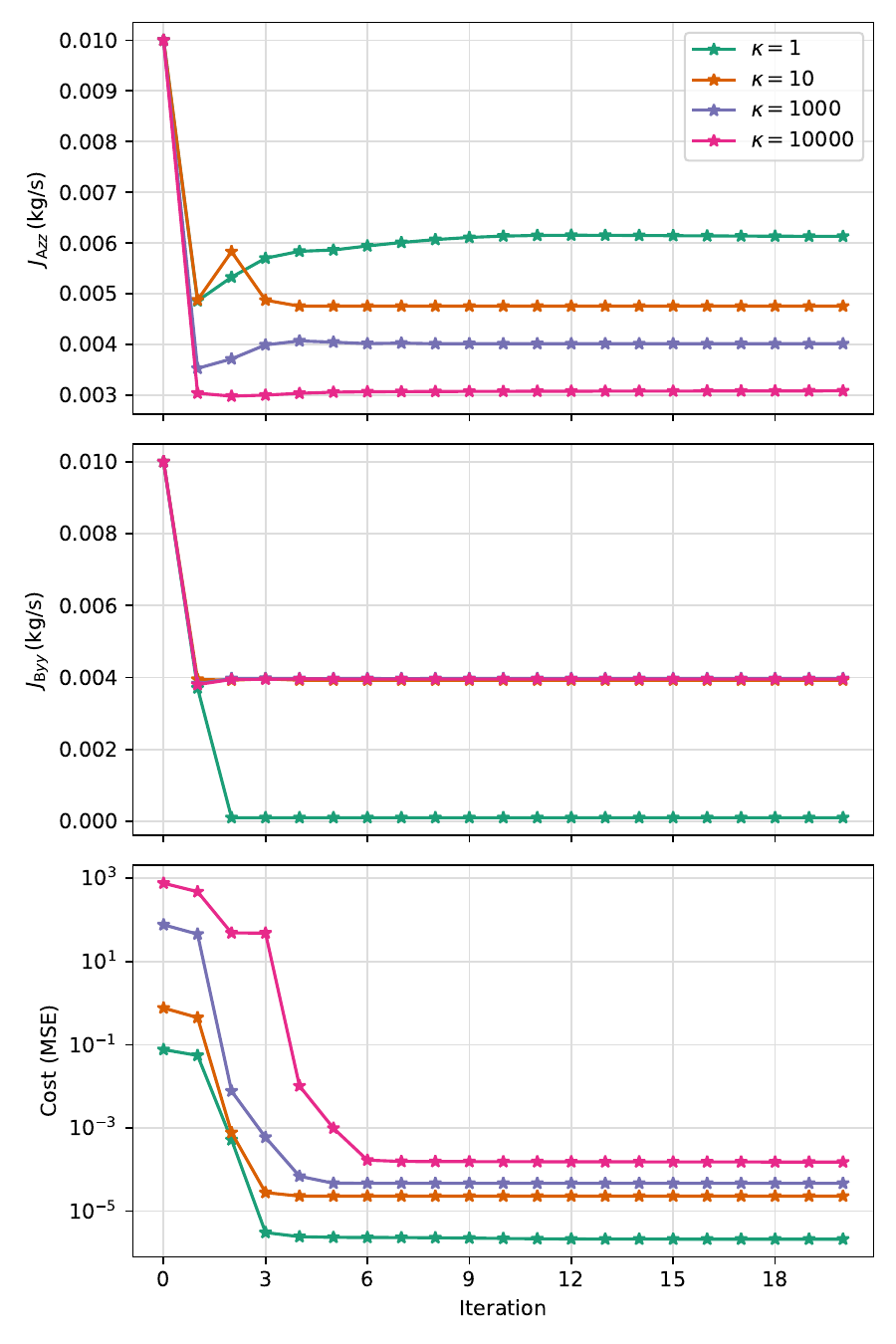}
    \caption{
        The left column of plots shows simulated observation trajectories using models tuned with different $\kappa$.
        From top to bottom, we have measured joint angles $\theta_1$, $\theta_2$, and the control signal $u$ versus time.
        The right column of plots shows the parameter values for two of the tuned parameters versus iteration for the same set of $\kappa$.
        The bottom right plot shows how the cost decreases under optimization.
    }
    \label{fig:furuta_rollout_weight_sweep}
\end{figure}

To put the errors into perspective, we evaluate the resulting models through simulation---we use them to reproduce the observation trajectory that they were tuned on.
We feed the simulation the recorded control signal.
We see four behaviors in \figref{fig:furuta_rollout_weight_sweep}:
The method overfits the observation data for small state error weights ($\stateweight = 1$), resulting in a poor simulation model.
For medium state error weights ($\stateweight = 10$), the optimization converges quickly and produces a model with good short-term simulation accuracy.
For large state error weights ($\stateweight = 1000$), the optimization converges slightly slower but produces a model with good long-term simulation accuracy.
The beginning observation deviates from the other trajectories for even larger state error weights ($\stateweight = 10000$).
The larger deviation in simulation error suggests that the optimization gets stuck in a local minima.

\figref{fig:furuta_rollout_weight_sweep} also shows how the cost and two of the inertia parameters change during optimization for different $\stateweight$.
The parameters converge quickly at the start, and, except for $\kappa=1$, the parameter values stabilize after roughly four iterations.
The cost is also reduced quickly, although it is a bit slower for higher $\kappa$ in this particular instance.
Still, not much happens after six iterations.
    
\begin{table}
\caption{
Estimated parameters of the real Furuta pendulum.
Values in parenthesis are the standard deviation estimates to the same precision as the least significant digit.
If the estimated parameter has the value $0$, then the values in parenthesis are the standard deviation in absolute precision.
}
\centering
\label{tab:furuta_baseline_mix}
\begin{tabular}{l|rrrrr|r}
\toprule
{} & {0} & {1} & {2} & {Combined} \\
\midrule
Cost & 3.22e-05 & 2.28e-06 & 3.01e-05 & 2.15e-05 \\
CV & 2.46e-05 & 2.36e-04 & 2.34e-05 & 9.47e-05 \\
$J_{\tA xx}\,(\si{\kilogram\meter^2})$ & 0(200.0) & 0(2000.0) & 0(200.0) &  \\
$J_{\tA zz}\,(\si{\kilogram\meter^2})$ & 0.00449(2) & 0.0042(2) & 0.00438(4) &  \\
$J_{\tA yy}\,(\si{\kilogram\meter^2})$ & 0(50.0) & 0(3000.0) & 0(80.0) &  \\
$J_{\tB yy}\,(\si{\kilogram\meter^2})$ & 0.003964(1) & 0.00384(2) & 0.003962(2) & \\
$J_{\tB zz}\,(\si{\kilogram\meter^2})$ & 0.00242(6) & 0.00411(2) & 0.00374(9) &  \\
$b_1\,(\si{\kilogram\per\second})$ & 0.0009(2) & 0.00140(6) & 0(0.0006) &  \\
$b_2\,(\si{\kilogram\per\second})$ & 0.00010(2) & 0(0.0001) & 0.00006(5) &  \\
$\mu_1$ & 0.00058(1) & 0(2e-05) & 0.00072(6) &  \\
$\mu_2$ & 0.00022(3) & 0.00009(2) & 0.00019(8) &  \\
$\kappa$ & - & 7.29(5) & 7.28(2) &  \\
\bottomrule
\end{tabular}
\end{table}

\subsection{Calibration with simulated data}
To further investigate the method's sensitivity, we perform parameter studies on synthetic observation data.
We generate this observation data by simulating the release scenario of a Furuta pendulum with timestep $h = 10^{-4}\,\si{s}$ and compliance $\epsilon = 10^{-6}\,\si{m/N}$.
We add zero mean Gaussian noise with standard deviation \SI{0.001}{rad} to the simulated joint angles to make the setting more realistic.
Further details of the settings used for simulation are found in \ref{appendix:data_sim_furuta_pendulum}.
We then use another simulation model with the same structure as the identification target.
For simplicity, only four parameters are subject to identification: $J_{\tA zz}$, $J_{\tB yy}$, $b_1$ and $b_2$.
We set the initial values of these parameters according to Table \ref{tab:furuta_default_parameters}.
The method parameters of interest are the state error weight $\stateweight$ and the timestep $h$.
We also investigate the sensitivity to the compliance $\epsilon$.

\begin{figure}
    \centering
    \includegraphics[width=0.495\textwidth]{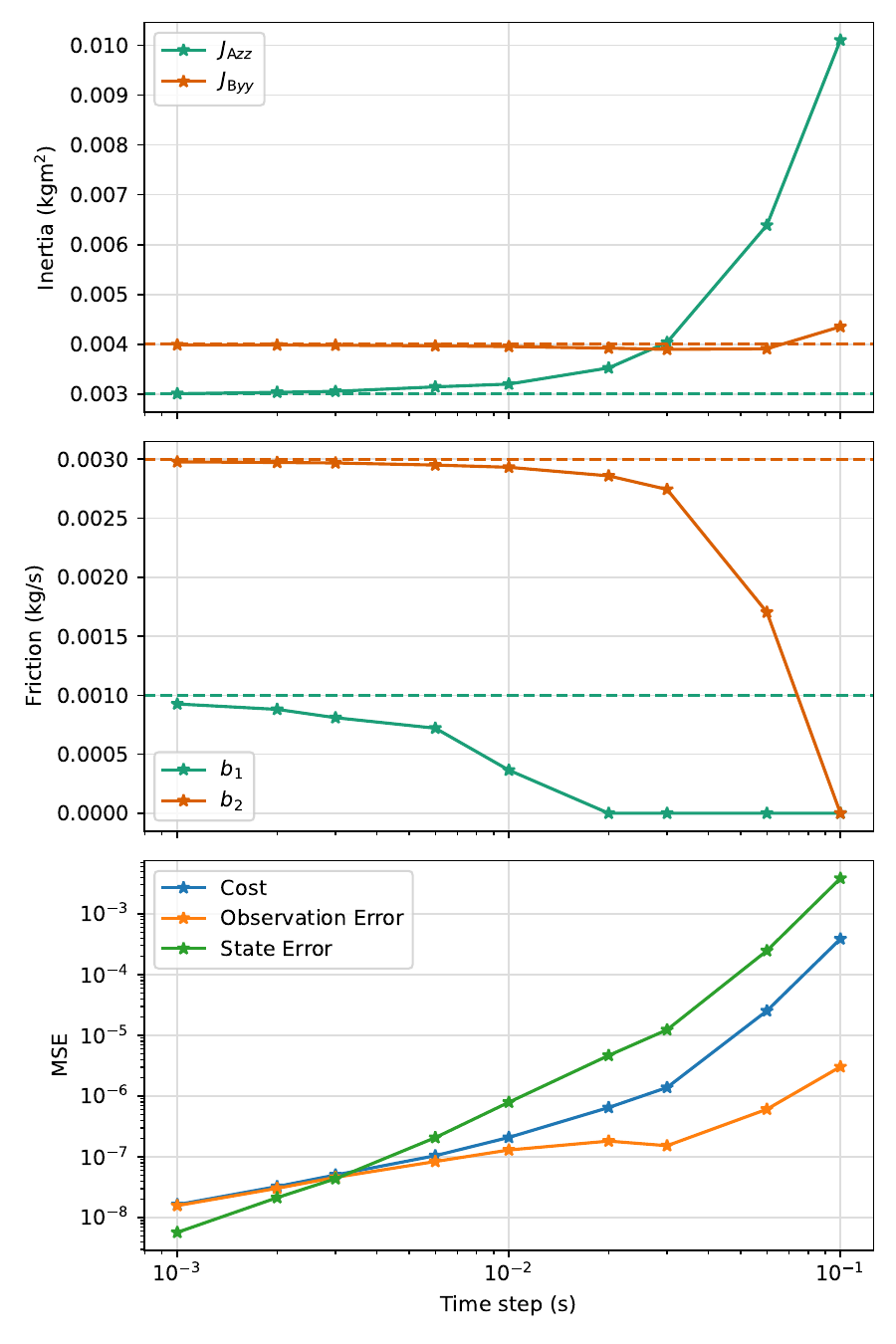}
    \includegraphics[width=0.495\textwidth]{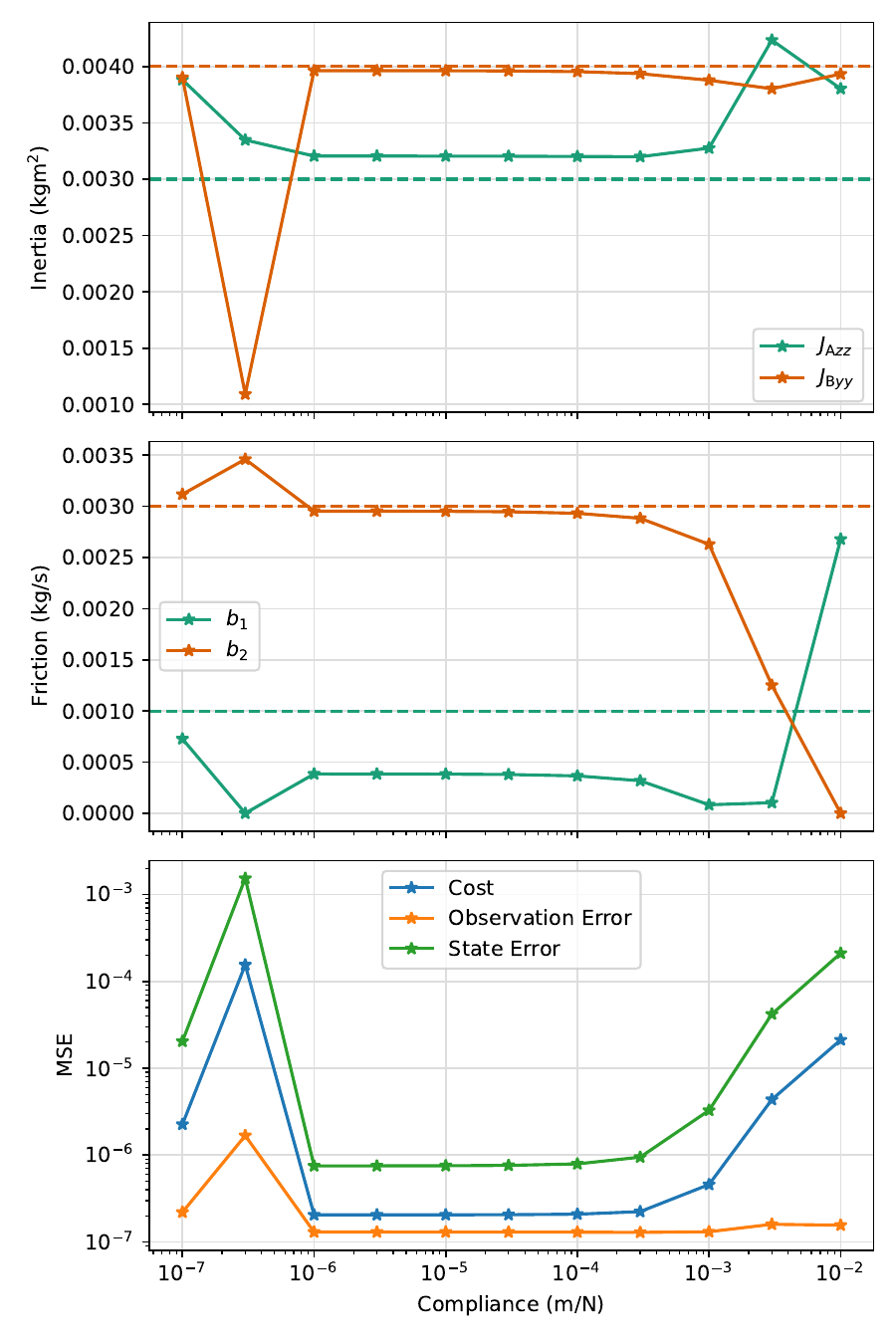}
    \caption{
        The left plot column shows the resulting parameters and errors after optimization with different time-steps.
        The right plot column is the same but plotted for different compliances.
    }
    \label{fig:compliance_timestep_sweeps}
\end{figure}

First, we identify the parameters for different choices of time-steps with fixed compliance, $\epsilon = 10^{-4}\,\si{m/N}$ and fixed state error weight, $\stateweight = 10^2$, see \figref{fig:compliance_timestep_sweeps}.
The identified parameters are found to weakly depend on the choice of time-step:
The parameter $b_1$ sees large relative changes, even going from $h = 10^{-2}\,\si{s}$ to $h = 10^{-3}\,\si{s}$.
For time-steps larger than \SI{2e-2}{s}, the estimated parameters change rapidly, and the errors increase further.

Second, the parameters are identified for different choices of compliance with fixed time-step, $h = 10^{-2}\,\si{s}$, and fixed state error weight, $\stateweight = 10^2$, see \figref{fig:compliance_timestep_sweeps}.
The story is much the same as in the calibration with real data.
The identified parameters have a weak dependence on the choice of compliance---as long as $\epsilon \lesssim 10^{-4}$.
As in the real scenario, strange behavior is seen for very stiff constraints, $\epsilon \lesssim 10^{-6}$.
Note that the estimated parameters do not approach the true values that were used to generate the observations.
This is explained by the choice of time-step, which is $h = 10^{-2}\,\si{s}$ for all models in the compliance parameter sweep.

\begin{figure}
    \centering
    \includegraphics[width=0.9\textwidth]{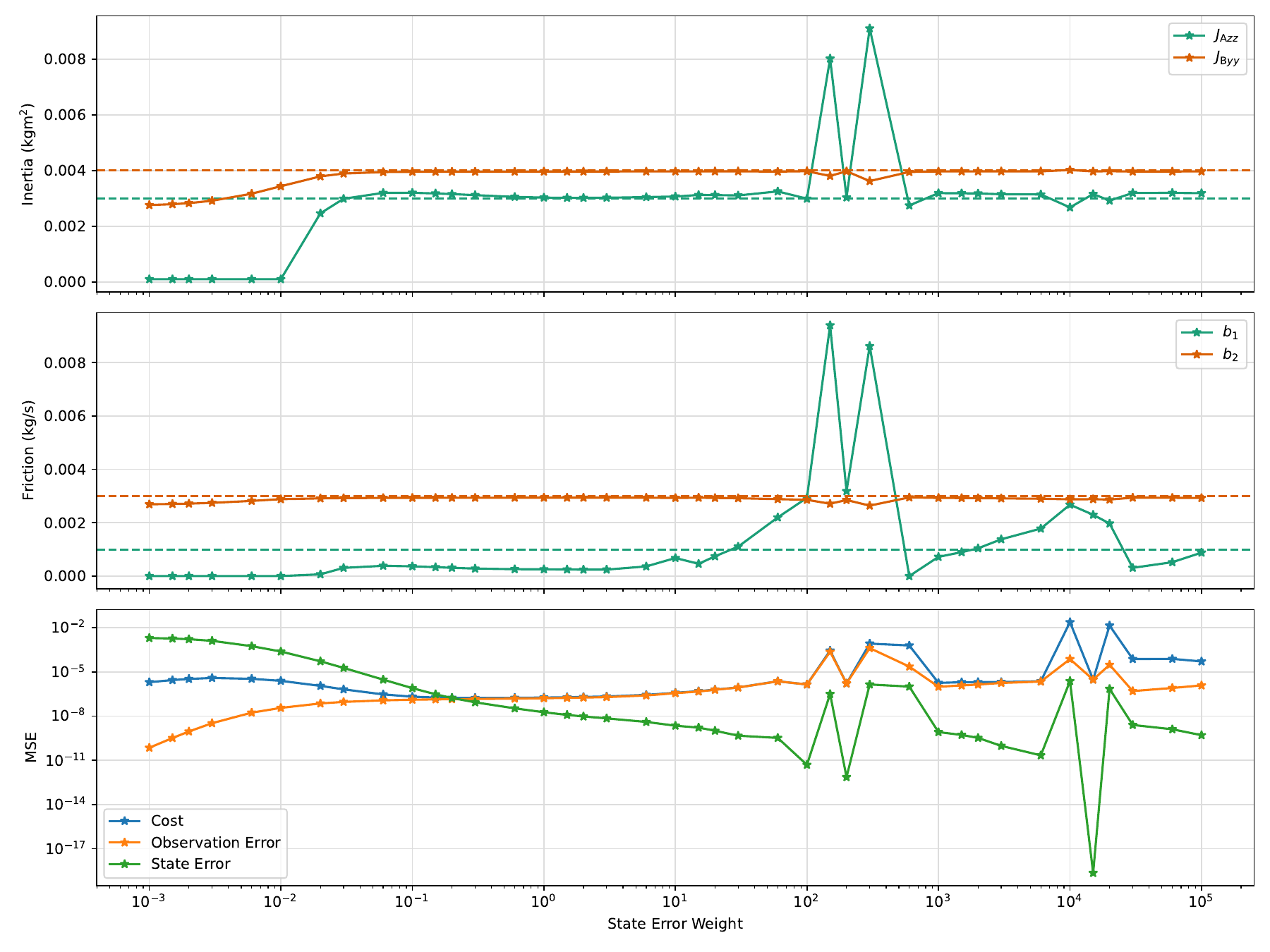}
    \caption{Paramer values versus and errors versus the state error weight in the simulated release scenario.}
    \label{fig:state_error_weight_sweep_sim}
\end{figure}

Finally, \figref{fig:state_error_weight_sweep_sim} shows the identified parameters for different choices of state error weight $\stateweight$.
In contrast to the real data, it is now enough to use a state error weight of $\kappa \approx 10^{-1}$ for the observation error to cross the simulation error and enter the plateau of stable parameters.
At $\kappa \approx 10^2$, the error starts to increase again.

\subsection{Computational time}
We use a desktop computer with an Intel Core i7-10700 processor at \SI{4.8}{GHz}.
In the example where we identify the parameters in Table \ref{tab:param_dry} for the swing-up scenario, it takes six iterations of Levenberg--Marquardt and \SI{6.6}{s} in real-time for the algorithm to converge.
This result is based on a single CPU core and should be taken as a reference value.
A more extensive study of computational time is outside the scope of this work.

\subsection{Discussion}
\label{sec:discussion}
\paragraph{Shadow estimation}
Parameter estimation in discrete dynamics is inherently different from that of a continuous model.
The estimated parameters generally differ from the physical parameters since discretization errors enter the optimization problem, as demonstrated in \figref{fig:compliance_timestep_sweeps}.
What we estimate can instead be considered \emph{shadow parameters} and \emph{shadow states} that fit the discretized model such that it \emph{shadows} the measurements.
We see this as a possible advantage---the parameters may be able to absorb some of the discretization errors that would otherwise be present.

\paragraph{Uncertainty estimation}
The inconsistencies of linearized uncertainty estimates over different time-series suggest that nonlinear effects are important.
In other words, more sophisticated uncertainty quantification is necessary to get reliable uncertainty estimates.
Still, linearized uncertainty gives an indication of the true uncertainty and is useful to determine non-excited parameters---parameters that makes no difference to the cost.

\paragraph{Inverse dynamics}
The results from the parameter study reveal ways in which our method may fail. Most of these failure modes are expected: If the time-step is too large, the dynamics is not captured. 
If the state error weight is too small, the model fits to noise. However, we also find that our method breaks down for sufficiently stiff constraints, see \figref{fig:state_error_weight_sweep_real} and \figref{fig:compliance_timestep_sweeps}.
We believe that this is well-explained by the way we compute the inverse dynamics. 
We solve for the inverse dynamics explicitly and not implicitly as we do for the constraint forces.
The reason to solve for the constraint forces implicitly in the first place is to achieve stability---dampen high-frequency dynamics that are of no interest to resolve.
The same logic could be applied to the impulse that corrects for errors.
In other words, we suspect that instabilities during optimization cause the derivatives to provide misdirection. An implicit method to do inverse dynamics could fix the problems at the cost of increased computational cost.

\section{Conclusions}
\label{sec:conclusions}
We develop a method for parameter estimation for regularized time-discrete multibody dynamics---a numerically robust formulation that allows deviations from the kinematic constraints.
These deviations result in unobserved state variables that appear in the discrete time-stepper.
To account for the unobserved state variables, we estimate both states and parameters by solving a nonlinear least squares optimization problem with a Levenberg--Marquardt algorithm.
We compute derivatives using forward mode accumulation, which our custom physics simulator supports.

We study the method's performance using synthetic and real measured data.
From the experiments, we conclude that joint state and parameter estimation through nonlinear least squares optimization and inverse dynamics can be a very efficient alternative to simulation error optimization.
The experimental results demonstrate, in particular, that this approach can identify the parameters in a 13-dimensional parameter space with impressive precision in just four iterations of Levenberg-Marquardt.

We identify two important method parameters:
First, the choice of state error weight.
If it is too large, the problem becomes infeasible; if it is too small, the method overfits the observations.
Second, the choice of compliance parameter for the kinematic constraints.
The optimizer will have difficulties making proper progress if the constraints are too stiff.
As a result, kinematic chains have to be sufficiently soft for the presented approach to work as expected.
Implicit formulations of inverse dynamics may address this problem, but this is outside the scope of the current work.

Future research directions include investigating the method from a statistical perspective to find ways of mitigating bias, investigating differentiability and optimization under non-smooth impacts, and reducing the dimensionality of the state residuals by enforcing the time-stepping scheme during optimization.
We would also like to extend numerical studies to complex systems with more components.

\section*{Acknowledgements}
This work was partially supported by LKAB (AC-438), eSSENCE, Kempe Foundation (SMK-2056, U56), and the Wallenberg AI, Autonomous Systems and Software Program (WASP) funded by the Knut and Alice Wallenberg Foundation.

\bibliographystyle{plain}
\bibliography{references}

\appendix
\section{MLCP Solver}
\label{appendix:solver}
Consider the following MLCP
\begin{equation}
    \label{eq:generic_mlcp}
    \begin{aligned}
        \underbrace{
        \begin{bmatrix}
            \mat{S}_* & \mat{S}_{\text{B}}^T     \\
            \mat{S}_{\text{B}} & \mat{S}_{\text{BB}} \\
        \end{bmatrix}
        }_\mat{S}
        \underbrace{
        \begin{bmatrix}
            \vec{\lambda}_\tA   \\
            \vec{\lambda}_\tB \\
        \end{bmatrix}}_{\vec{\lambda}}
         & +
        \underbrace{
        \begin{bmatrix}
            \vec{b}_\tA \\
            \vec{b}_\tB \\
        \end{bmatrix}}_{\vec{b}}
        =
        \underbrace{
        \begin{bmatrix}
            \vec{0}   \\
            \vec{w}_\tB \\
        \end{bmatrix}
        }_{\vec{w}}, \\
        \complementarityalign{\vec{\lambda}}{\vec{w}_\tB},
    \end{aligned}
\end{equation}
where $\mat{S}_*$ is a square symmetric, $\mat{S}_{\text{BB}}$ is square but not necessary symmetric and $\mat{S}_{\text{B}}$ has rows and columns of appropriate size.
To find a solution, we first consider the following subdivision of the matrix into smaller blocks:
\begin{equation}
    \label{eq:mlcp_schur_block}
    \begin{bmatrix}
        \mat{S}_* & \mat{S}_{\text{B}}^T     \\
        \mat{S}_{\text{B}} & \mat{S}_{\text{BB}} \\
    \end{bmatrix}
    =
    \begin{bmatrix}[ccc|c]
        \mat{S}_{11} & \hdots & \mat{S}_{1n} & \mat{S}_{1 \tB} \\
        \vdots       & \ddots & \vdots       & \vdots       \\
        \mat{S}_{n1} & \ddots & \mat{S}_{nn} & \mat{S}_{n \tB} \\
        \midrule
        \mat{S}_{\tB 1} & \hdots & \mat{S}_{\tB n} & \mat{S}_{\tB \tB} \\
    \end{bmatrix}
    =
    \mat{S}.
\end{equation}
To solve \eqref{eq:mlcp_schur_block}, we factor $\mat{S}$ as $\mat{S} = \mat{L}\mat{D}\mat{L}^T$, where
\begin{equation}
    \mat{L} = \begin{bmatrix}[ccc|c]
        \mat{I}      & \hdots & \mat{0}      & \mat{0} \\
        \vdots       & \ddots & \vdots       & \vdots  \\
        \mat{L}_{n1} & \ddots & \mat{I}      & \mat{0} \\
        \midrule
        \mat{L}_{\tB 1} & \hdots & \mat{L}_{\tB n} & \mat{I} \\
    \end{bmatrix}, \quad
    \mat{D} = \begin{bmatrix}[ccc|c]
        \mat{D}_{11} & \hdots & \mat{0}      & \mat{0}      \\
        \vdots       & \ddots & \vdots       & \vdots       \\
        \mat{0}      & \ddots & \mat{D}_{nn} & \mat{0}      \\
        \midrule
        \mat{0}      & \hdots & \mat{0}      & \mat{D}_{\tB \tB} \\
    \end{bmatrix},
\end{equation}
using block LDLT factorization---a generalization of Cholesky factorization.
This factorization is valid as long as the non-symmetric matrix $\mat{S}_{\text{BB}}$ is never used as a pivot.
Next, we solve $\mat{L}\vec{y} = -\vec{b} + \vec{w}$ for $\vec{y}$ where $\vec{y} = \mat{D}\mat{L}^T \vec{\lambda}$ using forward substitution.
Note that for $\vec{y}_\tB$, we have
\begin{equation}
    \vec{y}_B = \vec{w}_B - \vec{b}_B - \sum_{j=1}^{i-1} \vec{L_{ij}}\lambda_B.
\end{equation}
Since $\vec{w}_B$ is unknown at this point, we can only solve for $\vec{y}^* = \vec{y} - \vec{w}$ and the complementarity conditions are yet to be resolved.
The diagonal scaling phase amounts to
\begin{equation}
    \begin{aligned}
        \begin{bmatrix}
            \mat{D}_{*} & \mat{0}      \\
            \mat{0}      & \mat{D}_{\tB \tB}
        \end{bmatrix}
        \begin{bmatrix}
            \mat{L}_{*}^T & \mat{L}_{\tB}^T \\
            \mat{0}        & \mat{I}^T
        \end{bmatrix}
        \begin{bmatrix}
            \lambdaa \\
            \lambdab \\
        \end{bmatrix}
         & +
        \begin{bmatrix}
            \vec{y}_\tA \\
            \vec{y}_\tB \\
        \end{bmatrix}
        =
        \begin{bmatrix}
            \vec{0}   \\
            \vec{w}_B \\
        \end{bmatrix}, \\
        \complementarity{\lambdab}{\vec{w}_{\tB}}
    \end{aligned}.
\end{equation}
The problem now separates into a pure linear complementarity problem (LCP) and multiple small linear systems for the diagonal scaling phase.
This LCP is solved using a basic implementation of Lemke's algorithm, \cite{murty1988linear}.
After a backward substitution phase, the resulting multipliers are used in \eqref{eq:velocity_update} to get the velocity update.

\section{Jacobian of the residuals}
\label{appendix:jacobian}
The current goal is to use forward-mode accumulation to compute the Jacobian matrix of the residual function $\Delta \vec{z}$ with respect to the states $\vec{x}_{0:n-1}$ and parameters $\param$.
In theory, this is best done by computing the analytical derivatives.
Here, our philosophy is to use automatic differentiation to spare ourselves that extra work.
While not optimal from a performance perspective, the resulting Jacobian is still very usable in practice.

To exploit the known sparsity patterns, we first realize that a Jacobian-matrix product is just several instances of Jacobian-vector products neatly packed together.
Consider the upper block-bidiagonal Jacobian $\partial \Delta \vec{p} / \partial \vec{x}_{0:n-1} = \mat{P}$. The following Jacobian-matrix product reveals all non-zero elements of $\mat{P}$:
\begin{equation}
    \underbrace{
        \begin{bmatrix}
            \mat{P}_{11} & \mat{P}_{12} & \;     & \;     & \;     \\
            \;     & \mat{P}_{22} & \mat{P}_{23} & \;     & \;     \\
            \;     &        & \mat{P}_{33} & \mat{P}_{34} & \;     \\
            \;     & \;     & \;     & \ddots & \ddots \\
        \end{bmatrix}
    }_{\mat{P}}
    \underbrace{
        \begin{bmatrix}
            \mat{I}      & \mat{0}      \\
            \mat{0}      & \mat{I}      \\
            \mat{I}      & \mat{0}      \\
            \mat{0}      & \mat{I}      \\
            \vdots & \vdots \\
        \end{bmatrix}
    }_{\mat{E}}
    =
    \underbrace{
        \begin{bmatrix}
            \mat{P}_{11} & \mat{P}_{12} \\
            \mat{P}_{23} & \mat{P}_{22} \\
            \mat{P}_{33} & \mat{P}_{34} \\
            \vdots & \vdots
        \end{bmatrix}
    }_{\bar{\mat{P}}}
\end{equation}
That is, we can get the matrix $\bar{\mat{P}}$ and thus $\mat{P}$ by forward mode accumulation using Jacobian-vector products with columns from $\mat{E}$.
In fact, we can get $\mat{P}$ and $\mat{Y} = \partial \Delta \vec{y} / \partial \vec{x}_{0:n-1}$ (which is only block-diagonal) at the same time with
\begin{equation}
    \begin{bmatrix}
        \mat{Y} \\ \mat{P}
    \end{bmatrix}
    \mat{E}
    =
    \begin{bmatrix}
        \bar{\mat{Y}} \\ \bar{\mat{P}}
    \end{bmatrix}.
\end{equation}
To get the full Jacobian, we may use
\begin{equation}
    \begin{bmatrix}
        \mat{Y}_\chi & \mat{Y} \\
        \mat{P}_\chi & \mat{P} \\
    \end{bmatrix}
    \begin{bmatrix}
        \mat{I} & \mat{0} \\
        \mat{0} & \mat{E}
    \end{bmatrix}
    =
    \begin{bmatrix}
        \mat{Y}_\chi & \bar{\mat{Y}} \\
        \mat{P}_\chi & \bar{\mat{P}} \\
    \end{bmatrix},
\end{equation}
where $\mat{Y}_\chi = \partial \Delta \vec{y} / \partial \param$ and  $\mat{P}_\chi = \partial \Delta \vec{p} / \partial \param$.

There is one final problem to address with this approach.
The quaternions are represented by four scalars, while the tangent space of rotations only has three dimensions.
In other words, the above approach will not give the correct Jacobian for on-manifold optimization.
To address this, we consider new residual functions $\Delta \vec{p}_k'$ that are also functions of rotation vectors $\vec{\rv}_{\tA, 0:n-1}$, $\vec{\rv}_{\tB, 0:n-1}$, $\hdots$ for each rigid body and time-step.
These rotation vectors are additional increments to the rotations of all bodies.
By inputting $\vec{\rv}_{\tA, 0:n-1} = \vec{\rv}_{\tB, 0:n-1} = \hdots = \vec{0}$, these extra inputs can be used to get the three-dimensional tangent values without interfering with the value of the residuals.
The custom differentiation rule required to achieve this is explained in \ref{appendix:jvp_quaternion_from_rotation_vector}.

\section{Constrained Linear Least Squares}
\label{appendix:box_constrainted_lls}
In each iteration of the Levenberg-Marquardt algorithm, we need to solve a constrained linear least squares problem of the type
\begin{mini}|l|
    {\vec{\delta} \bar{\param}}{
        \norm{\combresiduals + \mat{J} \delta \bar{\param}}^2_{\mat{W}} + \mu \|\delta \bar{\param}\|^2
    }
    {\label{eq:c_linear_least_squares}}{}
    \addConstraint{\param + \delta \param \in B}.
\end{mini}
We restrict $B$ to be a box ($B$ as a convex polytope is still applicable), making us able to write the constraint as $\mat{A}\vec{\delta} \bar{\param} \geq \vec{b}$.
Now, introduce a slack variable $\vec{w} \geq \vec{0}$ such that $\mat{A} \vec{\delta} \bar{\param} = \vec{h} + \vec{w}$.
The Lagrangian is
\begin{equation}
    \mathcal{L}(\delta \bar{\param}, \vec{\lambda}) = \norm{ \mat{J}\delta \bar{\param} + \vec{b}}^2_{\mat{W}} + \mu \|\delta \bar{\param}\|^2 + \vec{\lambda}^T (\mat{A} \delta \bar{\param} - \vec{b} - \vec{w}).
\end{equation}
By the Karush-Kuhn-Tucker conditions, we can restate the optimality conditions as an MLCP
\begin{equation}
    \label{eq:bc_linear_least_squares_mlcp}
    \begin{aligned}
        \begin{bmatrix}
            \mat{J}^T\mat{W}\mat{J}+\mu \mat{I} & \mat{A}^T     \\
            \mat{A} & \mat{0} \\
        \end{bmatrix}
        \begin{bmatrix}
            \delta \bar{\param} \\
            \vec{\lambda} \\
        \end{bmatrix}
         & +
        \begin{bmatrix}
            \mat{A}^T \mat{W} \vec{b} \\
            -\vec{b} \\
        \end{bmatrix}
        =
        \begin{bmatrix}
            \vec{0} \\
            \vec{w} \\
        \end{bmatrix}, \\
        \complementarityalign{\vec{\lambda}}{\vec{w}}.
    \end{aligned}
\end{equation}
This MLCP is solved using the method outlined in \ref{appendix:solver}.

\section{Custom differentation rules}
\label{appendix:jvp_rules}
In this section, we use dot notation to denote tangent vectors or Jacobian vector products (JVPs):
\begin{equation}
    \dot{\vec{w}} = \frac{d \vec{w}}{d x},
\end{equation}
and should not be confused with time derivatives.
\subsection{JVP of the linear complementarity problem}
Consider the following LCP:
\begin{equation}
    \label{eq:lcp}
    \begin{aligned}
        \mat{H} \vec{z} + \vec{q} &= \vec{w}, \\
        \complementarityalign{\vec{w}}{\vec{q}}.
    \end{aligned}
\end{equation}
We treat the solver to this LCP can be treated as a function: $\vec{z} = \text{LCP}(\mat{H}, \mat{q})$.
When augmented with tangents, the solver will, in addition to $\mat{H}$ and $\vec{q}$, be passed the tangents $\dot{\mat{H}}$ and $\dot{\vec{q}}$, and should return the tangent $\dot{\vec{z}}$ in addition to $\vec{z}$.
By the complementarity conditions, we can reformulate \eqref{eq:lcp} terms of active indices $\alpha$ and passive indices $\beta$:
\begin{equation}
    \label{eq:lcp_in_complementary_basis}
    \mat{A} \vec{y} =
    \begin{bmatrix}
        -\mat{H}_{\bullet \alpha} & \mat{I}_{\bullet \beta}
    \end{bmatrix}
    \begin{bmatrix}
        \vec{z}_\alpha \\ \vec{w}_\beta
    \end{bmatrix}
    =
    \vec{q}.
\end{equation}
Implicit differentiation yields
\begin{equation}
    \dot{\mat{A}}\vec{y} + \mat{A}\dot{\vec{y}} = \dot{\vec{y}}.
\end{equation}
We can rewrite $\dot{\mat{A}} \vec{y}$ as
\begin{equation}
    \dot{\mat{A}}\vec{y} =     \begin{bmatrix}
        -\dot{\mat{H}}_{\bullet \alpha} & \mat{0}_{\bullet \beta}
    \end{bmatrix}\vec{y}=
    -\dot{\mat{H}}_{\bullet \alpha} \vec{z}_\alpha = -\dot{\mat{H}}  \vec{z},
\end{equation}
to have a linear system for $\dot{\vec{y}}$ with a right-hand side of known quantities
\begin{equation}
   \mat{A} \dot{\vec{y}} = \dot{\vec{q}} + \dot{\mat{H}} \vec{z}.
\end{equation}
The solution $\vec{z}$ is related to $\vec{y}$ by
\begin{equation}
    \vec{z} =
    \begin{bmatrix}
        -\vec{z}_\alpha \\ 0
    \end{bmatrix}
    =
    \begin{bmatrix}
        \mat{I}_{\alpha\alpha} & \mat{0} \\
        \mat{0}                & \mat{0}
    \end{bmatrix}
    \begin{bmatrix}
        \vec{z}_\alpha \\ \vec{w}_\beta
    \end{bmatrix}
    =
    \vec{1}_\alpha \odot \vec{y},
\end{equation}
where $\vec{1}_\alpha$ is the vector with ones at indices $\alpha$ and zeros elsewhere, and $\odot$ is the element-wise vector product.
The tangent is $\dot{\vec{z}} = \vec{1}_\alpha \odot \dot{\vec{y}}$.

\subsection{JVP of quaternion from rotation vector}
\label{appendix:jvp_quaternion_from_rotation_vector}
\def\rv{\psi}
Zero-valued rotation vectors $\vec{\psi}=\vec{0}$ are used to carry tangent information of rotations during optimization. 
This enables the construction of the correct Jacobian as well as updates in the form of rotation vectors $\vec{\delta}_{\vec{\psi}}$.
These rotation vectors need to be converted to quaternions, creating a need for custom differentiation rules for this operation.

Consider a quaternion $\vec{e}$ that is constructed from a rotation vector $\vec{\rv} = \hat{\vec{n}}\rv$, where $\hat{\vec{n}}$ is a unit vector.
We write
\begin{equation}
    \vec{\quat} = \begin{bmatrix} \quat_s \\ \vec{\quat}_v \end{bmatrix} = \begin{bmatrix} \cos{\frac{\rv}{2}} \\ \sin{\frac{\rv}{2}} \hat{\vec{n}} \end{bmatrix}.
\end{equation}
The tangent of the scalar component $\quat_s$ is
\begin{equation}
    \dot{\quat}_s = \frac{\vec{\rv} \cdot \dot{\vec{u}}}{2 \rv} \sin{\frac{\rv}{2}}
\end{equation}
and the tangent of the vector component $\quat_v$ is
\begin{equation}
    \dot{\vec{\quat}}_v = \Big(\frac{ \dot{\vec{\rv}}}{\rv} - \frac{\vec{u} \cdot \dot{\vec{\rv}}}{u^3}\Big)\sin{\frac{u}{2}} +  \hat{n} \frac{\vec{u} \cdot \dot{\vec{\rv}}}{2\rv} \cos{\frac{u}{2}}.
\end{equation}
To avoid division by zero, we look at the series expansion of $\sin$ and $\cos$:
\begin{equation}
    \begin{aligned}
        \dot{\vec{\quat}}_v &= \Big(\frac{ \dot{\vec{\rv}}}{\rv} - \frac{\vec{\rv} \cdot \dot{\vec{\rv}}}{\rv^3}\Big)\Big(\frac{\rv}{2} - \frac{\rv^3}{2^3 3!} +\mathcal{O}(\rv^5)\Big) +  \hat{n} \frac{\vec{\rv} \cdot \dot{\vec{\rv}}}{2\rv} \Big( 1 - \frac{\rv^2}{2^2 2!} + \mathcal{O}(\rv^4)\Big)
        \\&= \frac{\dot{\vec{\rv}}}{2} - \frac{1}{24} \vec{\rv}(\vec{\rv} \cdot \dot{\vec{\rv}})  + \mathcal{O}(\rv^2), \\
        \dot{\vec{\quat}}_s &= \frac{\vec{\rv} \cdot \dot{\vec{\rv}}}{2 \rv}(\frac{\rv}{2} - \frac{\rv^3}{2^3 3!} +\mathcal{O}(\rv^5)\Big) = \vec{\rv} \cdot \dot{\vec{\rv}} \Big( \frac{1}{4} - \frac{\rv^2}{2^4 3!} - \mathcal{O}(\rv^4)\Big).
    \end{aligned}
\end{equation}
If we know that $\vec{\rv}=\vec{0}$, the expression reduces to
\begin{equation}
    \dot{\vec{\quat}} = \begin{bmatrix} 0 \\ \frac{1}{2} \dot{\vec{\rv}} \end{bmatrix}.
\end{equation}

\section{Maximum likelihood}
\label{appendix:maximum_likelihood}
The maximum likelihood estimation method for estimating the state variables and parameters amounts to finding
\begin{equation}
    \hat{\vec{x}}_{1:n}, \hat{\param} = \argmax_{\vec{x}_{1:n}, \param} p(\vec{x}_{1:n}, \param | \vec{y}_{1:n}).
\end{equation}
We rewrite the likelihood function using Bayes' theorem for probability distributions
\begin{equation}
    \label{eq:bayes}
    p(\vec{x}_{1:n}, \param | \vec{y}_{1:n}) = \frac{p(\vec{y}_{1:n} | \vec{x}_{1:n}, \param) p(\vec{x}_{1:n}, \param)}{p(\vec{y}_{1:n})}.
\end{equation}
Note that $p(\vec{y}_{1:n})$ is independent of $\vec{x}_{1:n}$ and $\param$, so it can be dropped without changing the objective:
\begin{equation}
    \hat{\vec{x}}_{1:n}, \hat{\param} = \argmax_{\vec{x}_{1:n}, \param} p(\vec{y}_{1:n} | \vec{x}_{1:n}, \param) p(\vec{x}_{1:n}, \param).
\end{equation}
As the observation at time $k$ only depends on the state at time $k$, we have
\begin{equation}
    \label{eq:observations_prob_dist}
    \begin{aligned}
        p(\vec{y}_{1:n} | \vec{x}_{1:n}, \param) & =  \prod_{k=1}^n p(\vec{y}_k| \vec{x}_k, \param).
    \end{aligned}
\end{equation}
The joint probability distribution can rewritten as
\begin{equation}
    p(\vec{x}_{1:n}, \param) = p(\vec{x}_{1:n} | \param) p(\param).
\end{equation}
Now, since the state at time $k$ only depends on the state at time $k-1$, we have
\begin{equation}
    \begin{aligned}
        p(\vec{x}_{1:n} | \param) & = p(\vec{x}_n | \vec{x}_{1:n-1}, \param) p(\vec{x}_{1:n-1}, \param)       \\
                            & = p(\vec{x}_n | \vec{x}_{n-1}, \param) p (\vec{x}_{1:n-1}, \param)        \\
                            & = p(\vec{x}_1 | \param)\prod_{i=2}^n p(\vec{x}_i| \vec{x}_{i-1}, \param).
    \end{aligned}
\end{equation}
We will assume that $p(\vec{x}_1 | \param)$ is a known a prioi.
Substituting back in \eqref{eq:bayes}, we have
\begin{equation}
    p(\vec{x}_{1:n} | \vec{y}_{1:n}) = c p(\vec{x}_1 | \param) p(\param)\prod_{k=2}^n p(\vec{x}_k| \vec{x}_{k-1}, \param) \prod_{k=1}^n p(\vec{y}_k| \vec{x}_k, \param),
\end{equation}
where $c$ comes from the denomenator $p(\vec{y}_{1:n})$ of \eqref{eq:bayes}.
By taking the logarithm, we arrive at
\begin{align}
    -&\log{p(\vec{x}_{1:n} | \vec{y}_{1:n})} = c + \log{p(\param)} + \log{p_{x_0}(\vec{x}^*_1 - \vec{x}_1)} \nonumber \\
    &+ \sum_{k=2}^n \log{p_{\eta_k}(\vec{f}(\vec{x}_{k-1}) - \vec{x}_{k})}
    + \sum_{k=1}^n \log{p_{\xi_k}(\vec{h}(\vec{x}_k) - \vec{y}_{k})}.
\end{align}
If the distributions $p_{\eta_k}$ and $p_{\xi_k}$ for $k=1,\hdots,n$ are zero-mean Gaussian
\begin{equation}
    \begin{aligned}
        p_{\eta_k}(\vec{x}) = \frac{1}{\sqrt{(2\pi)^{n} |\mat{Q}|}} \exp \Big({-\frac{1}{2} \vec{x}^T \mat{Q}^{-1}_k \vec{x}} \Big), \\
        p_{\xi_k}(\vec{x}) = \frac{1}{\sqrt{(2\pi)^{n} |\mat{R}|}} \exp \Big({-\frac{1}{2} \vec{x}^T \mat{R}^{-1}_k \vec{x}} \Big),
    \end{aligned}
\end{equation}
then the non-constant terms are proportional to
\begin{equation}
    J(\vec{x}_{1:n}) =  \|\vec{x}^*_1 - \vec{x}_1\|^2_{\mat{Q}_{1}^{-1}}  + \sum_{k=2}^n \|\vec{f}(\vec{x}_{k-1}) - \vec{x}_k\|^2_{\mat{Q}_i^{-1}} + \sum_{k=1}^n \|(\vec{h}(\vec{x}_i) - \vec{y}_i)\|^2_{\mat{R}_i^{-1}} + \rho(\param),
\end{equation}
where $\rho$ is some penalty on $\param$ takes its form depending on $p(\param)$.
That is, maximizing the likelihood function with uniform $p(\param)$ is equivalent to solving a non-linear least squares objective.

\section{Frictional hinge joint}
\label{appendix:hinge}
The frictional hinge joint is a package of several kinematic constraints.
The basic parameters are compliance, damping, and location and orientation of the two attachment frames.
The frictional hinge joint optionally requires a viscous drag coefficient $b$, a dry frictional constant $\mu$ associated with the internal geometry as captured a distance $r$, and motor shaft inertia $J_{\text{shaft}}$.

The starting point is a holonomic hinge constraint.
Consider two coordinate systems, $\text{a}$ and $\text{b}$, attached to rigid bodies or the global frame. We impose the following restriction
\begin{equation}
    \label{eq:holonomic_hinge}
    \vec{0} = \vec{g}(\vec{q}) =
    \begin{bmatrix}
        \vec{r}_a - \vec{r}_b     \\
        \vec{e}_x^{(a)} \cdot \vec{e}_y^{(b)} \\
        \vec{e}_x^{(a)} \cdot \vec{e}_z^{(b)}
    \end{bmatrix},
\end{equation}
where the rotation matrices $\mat{R}_a = \begin{bmatrix}\vec{e}_x^{(a)} & \vec{e}_y^{(a)} & \vec{e}_z^{(a)}\end{bmatrix}$ and $\mat{R}_b = \begin{bmatrix}\vec{e}_x^{(b)} & \vec{e}_y^{(b)} & \vec{e}_z^{(b)}\end{bmatrix}$ describe the orientation of the frames $a$ and $b$.
Note that the compliance parameter will have different units when it is connected to translational and rotational constraints.
As such, a frictional hinge joint has two compliance parameters.

We use a regularized non-holonomic constraint to model velocity proportional friction (viscous friction).
The viscous friction coefficient is labeled $b$.
The friction should act in the rotational plane spanned by the null space of $\mat{G}$ (the Jacobian of $\vec{g}$).
To this end, we let $\bar{\mat{G}} = \begin{bmatrix} \vec{0}_{1 \times 3} & (\vec{e}_x^{(b)})^T\end{bmatrix}$.
The equation we impose to achieve viscous drag is
\begin{equation}
    \bar{\mat{G}}_k \vec{v}_{k+1} + \frac{1}{b h}\lambda_b = 0.
\end{equation}
We also include an option to control the hinge by a motor, modeled as an additional non-holonomic constraint.
To account for extra inertia of the shaft, $J_{\text{shaft}}$, we model the motor by
\begin{eqnarray}
    J_{\text{shaft}}^{-1} \lambda_{\text{motor}} + \bar{\mat{G}} \vec{v}_{k+1}= \bar{\mat{G}} \vec{v}_{k} + J_{\text{shaft}}^{-1} \tau_{\text{motor}},
\end{eqnarray}
where $\lambda_{\text{motor}}$ is a Lagrange multiplier.
To model Coulomb friction, we replace the holonomic constraint \eqref{eq:holonomic_hinge} with two contact constraints based on the same function with opposite sign, $\vec{g}(\vec{q}) \geq 0$ and $-\vec{g}(\vec{q}) \geq 0$, and construct the MLCP model described in Section \ref{sec:complementarity}.
We form $\bar{\mat{D}}$ by $\bar{\mat{D}} = \begin{bmatrix}
    \mat{G}^T & -\mat{G}^T
\end{bmatrix}^T$ and the row vector of friction coefficients could either have unique entries or be grouped in the following way:
\begin{eqnarray}
    \mat{U} = \begin{bmatrix}
        r \mu & r \mu & r \mu & \mu & \mu,
    \end{bmatrix}
\end{eqnarray}
where $r$ is a distance parameter that depends on the internal geometry of the joint.

\section{Data}
\subsection{Pendulum}
\label{appendix:pendulum}
\begin{figure}[ht]
    \includegraphics[width=\linewidth]{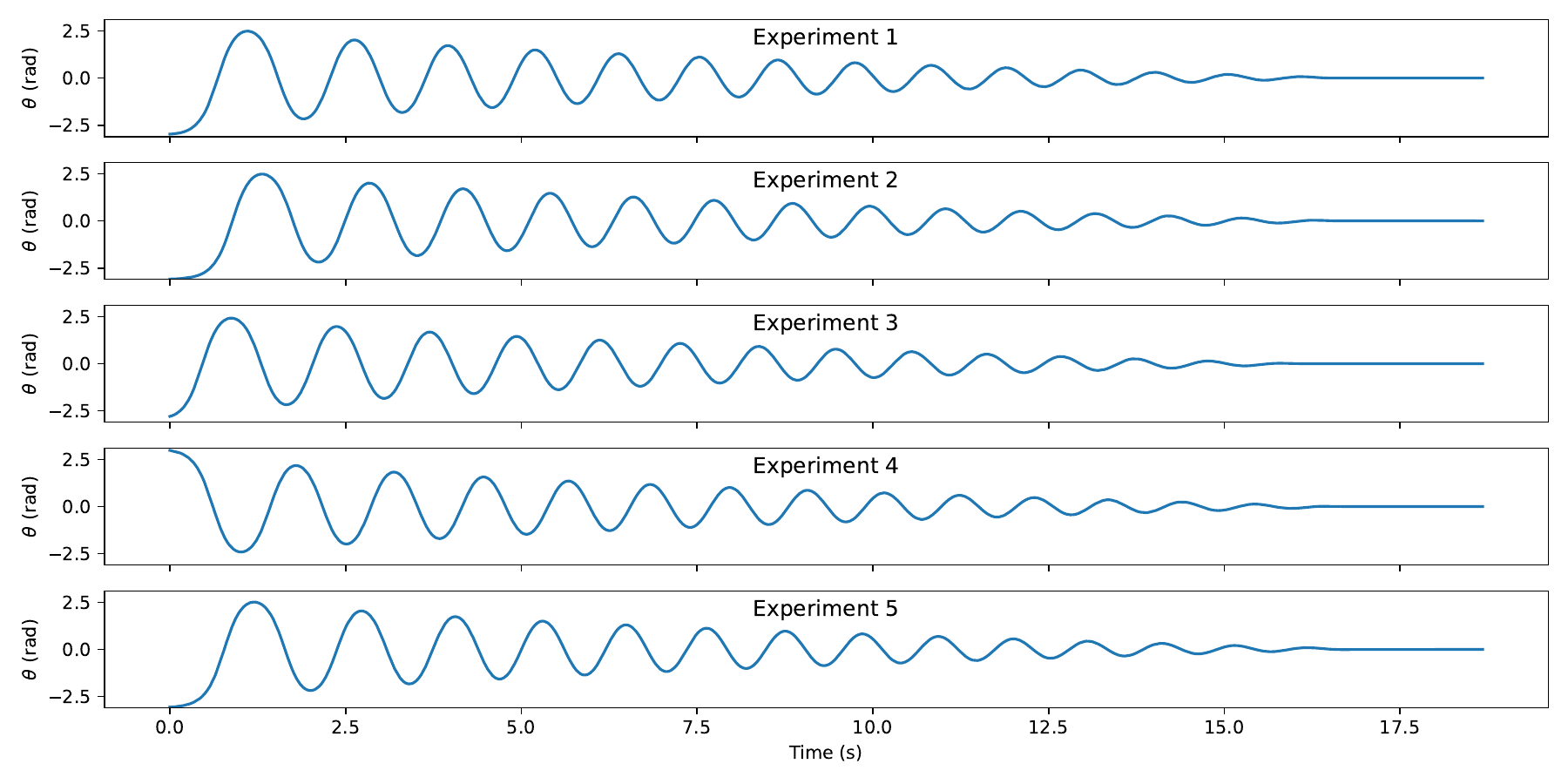}
    \caption{Observation data for all five time-series of the real pendulum.}
    \label{fig:pendulum_time_series}
\end{figure}

All five pendulum time-series are recorded by releasing the pendulum from an almost upright position.
The joint angles over time can be seen in \figref{fig:pendulum_time_series}.

\subsection{Real Furuta pendulum}
\label{appendix:data_real_furuta_pendulum}
\begin{figure}[ht]
    \centering
    \includegraphics[width=0.7\linewidth]{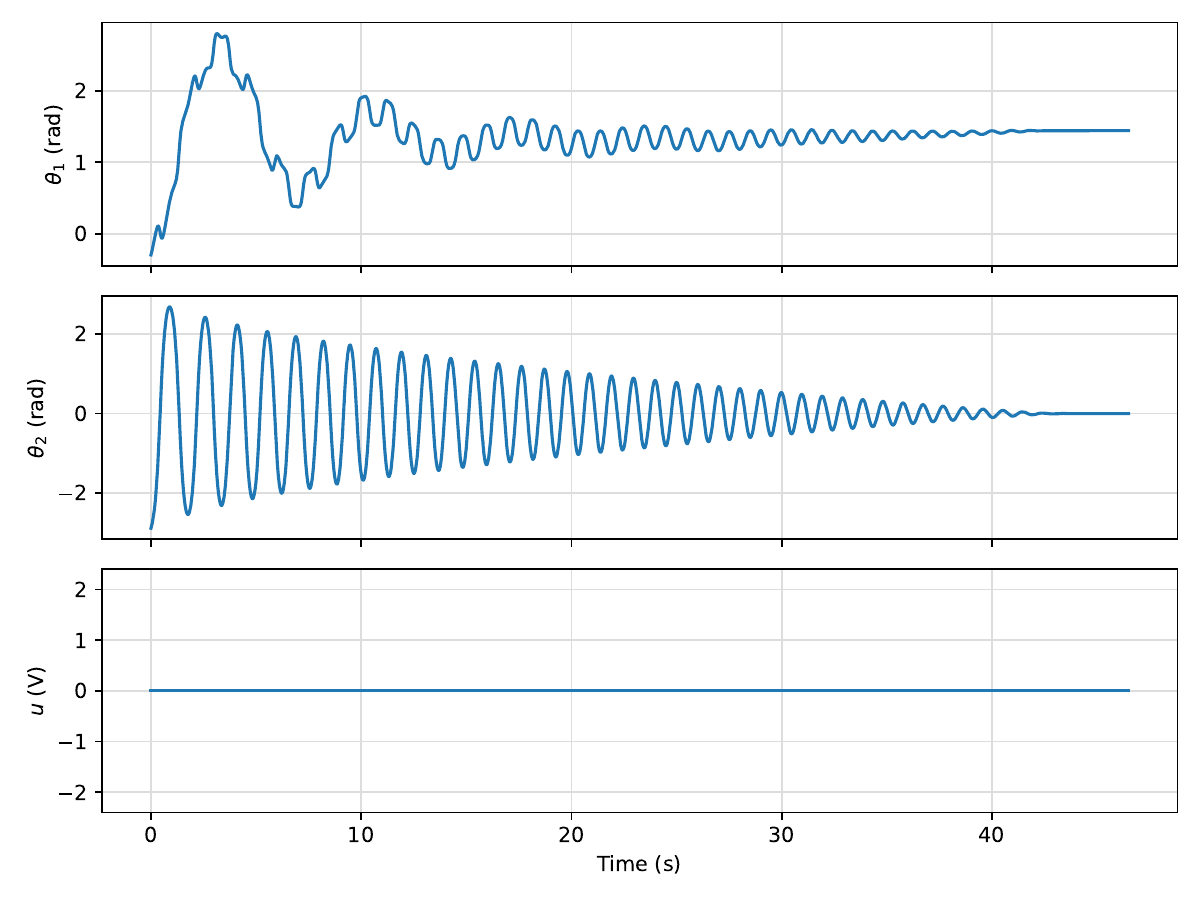}
    \caption{Observation data and control signal of the real Furuta pendulum release scenario.}
    \label{fig:furuta_release}
\end{figure}

\begin{figure}[ht]
    \centering
    \includegraphics[width=0.7\linewidth]{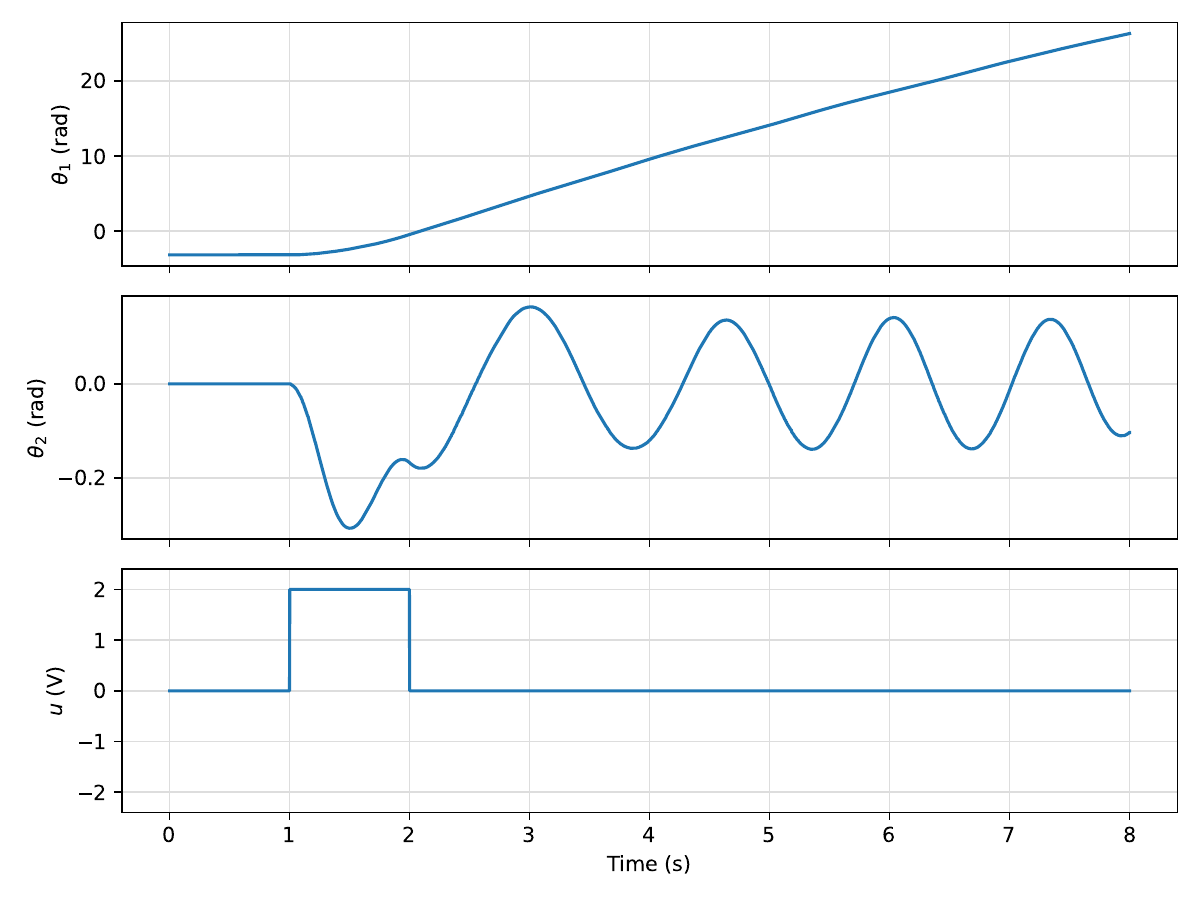}
    \caption{Observation data and control signal of the real Furuta pendulum pulse scenario.}
    \label{fig:furuta_pulse}
\end{figure}

Three different time-series of observations are collected using the real Furuta pendulum.
The time-series for the release scenario is recorded by releasing the Furuta pendulum from an upright position.
The data can be seen in \figref{fig:furuta_release}.
The time-series for the pulse scenario is recorded by using a one-second pulse as the control signal.
The resulting data and control signal is seen in \figref{fig:furuta_pulse}.

Finally, the swing-up scenario is recorded by using a closed-loop energy-shaping controller.
The resulting data and control signal is seen in \figref{fig:furuta_pulse}.
The details of the energy-shaping controller are left out as it is the recorded control signal, as seen in \figref{fig:furuta_pulse}, that is used for open-loop control in the parameter estimation experiments.

\begin{figure}[ht]
    \centering
    \includegraphics[width=0.7\linewidth]{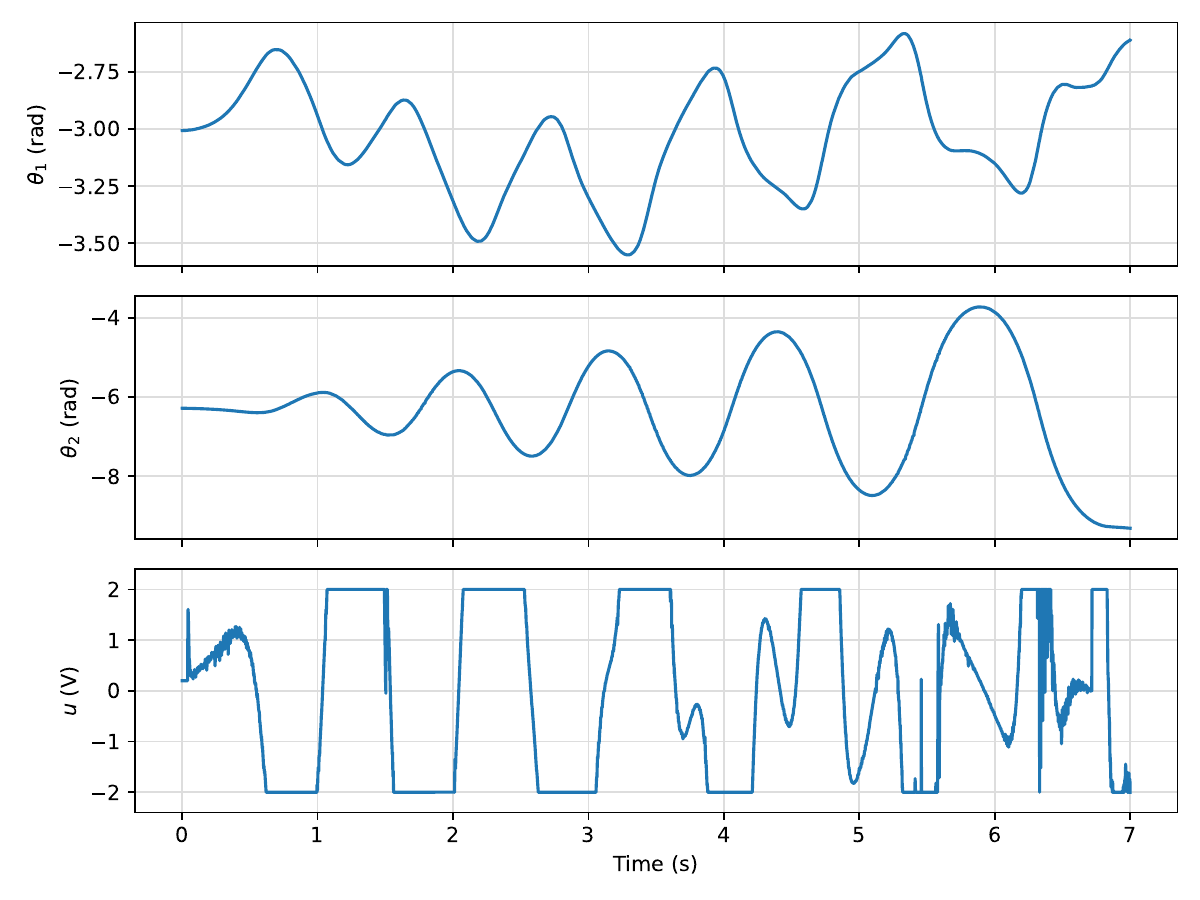}
    \caption{Observation data and control signal of the real Furuta pendulum swing-up scenario.}
    \label{fig:furuta_swingup}
\end{figure}

\subsection{Simulated Furuta pendulum}
\label{appendix:data_sim_furuta_pendulum}
\begin{figure}[ht]
    \centering
    \includegraphics[width=0.7\linewidth]{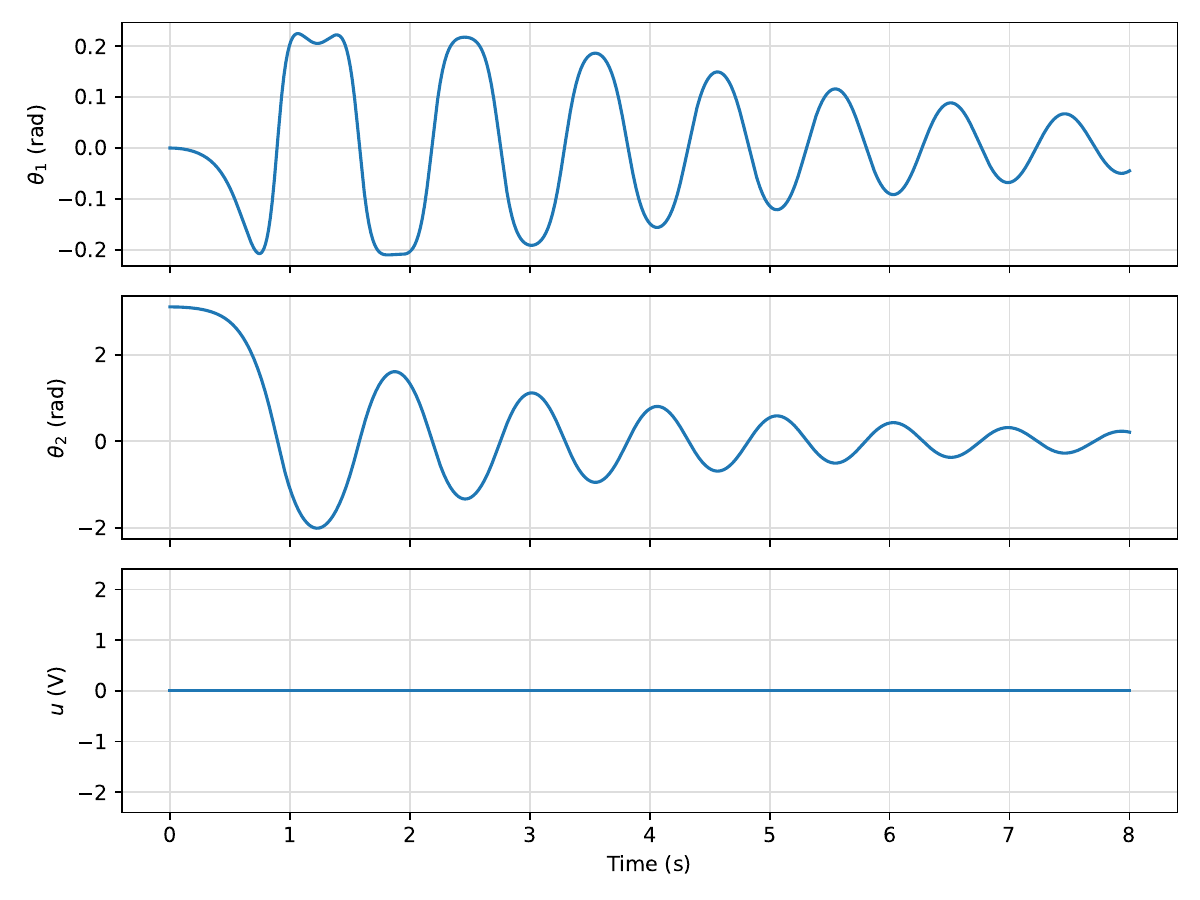}
    \caption{Synthetic observation data from the simulated Furuta pendulum.}
    \label{fig:furuta_realease_sim}
\end{figure}

The Furuta pendulum model that we use to generate the synthetic data has the parameters seen in Table \ref{tab:furuta_sim_parameters}.
The time-step is set to $h = 10^{-4}\,\si{s}$ and the initial configuration is at $\theta_1 = 0$ and $\theta_2 = (\pi - 0.03) \,\si{rad}$ at zero constraint violation and zero velocity.
The simulated trajectory is seen in \figref{fig:furuta_realease_sim}.

\begin{table}
\caption{Parameters used to generate the synthetic measurement data.}
\label{tab:furuta_sim_parameters}
\centering
\begin{tabular}{l|rrrrr|r}
    \midrule
    $l_1$ & \SI{0.128}{m} \\
    $l_\tA$ & \SI{0.248}{m} \\
    $l_2$ & \SI{0.92}{m} \\
    $\tau_1, \tau_2$ & $0.02\,\si{s}$ \\
    $\epsilon_1^{(\ell)}, \epsilon_2^{(\ell)}$ & $10^{-6}\,\si{m/N}$ \\
    $\epsilon_1^{(r)}, \epsilon_2^{(r)}$ & $10^{-6}\,\si{rad/N}$ \\
    $b_{1}$ & $10^{-3}\,\si{kg s^{-1}}$ \\
    \bottomrule
    \end{tabular}
    \begin{tabular}{l|rrrrr|r}
        \midrule
        $b_{2}$ & $10^{-4}\,\si{kg s^{-1}}$ \\
        $m_\tA$ & \SI{0.238}{kg} \\
        $m_\tB$ & \SI{0.428}{kg} \\
        $J_{\text{shaft}}$ & \SI{1e-4}{kg m^2} \\
        $\diag{(\mat{J}_{A})}$ & $\begin{bmatrix}
            10^{-4} & 0.003 & 0.003
        \end{bmatrix}^T\,\si{kg m^2}$ \\
        $\diag{(\mat{J}_{B})}$ & $\begin{bmatrix}
            10^{-4} & 0.004 & 0.004
        \end{bmatrix}^T\,\si{kg m^2}$ \\
        $\vec{g}_{\text{acc}}$ & $\begin{bmatrix}
            0 & -9.82 & 0
        \end{bmatrix}^T$\,\si{ms^{-2}} \\
        \bottomrule
        \end{tabular}
\end{table}

\end{document}